\pdfoutput=1
\documentclass[]{amsart}
\usepackage[T1,T2A]{fontenc}
\usepackage[utf8]{inputenc}
\usepackage{amsmath,mathrsfs,eucal}
\usepackage{mathtools} 
\usepackage{xcolor}
\usepackage{bm}
\usepackage{enumitem}
\usepackage{graphics} 
\usepackage{epsfig} 
\usepackage{graphicx,subfig}                     
\usepackage{wrapfig}
\usepackage{textcomp}
\usepackage{tikz}
\usepackage{pgfplots}
\usepackage{dsfont}
\pgfplotsset{compat=1.16}
\usepackage[cal=boondox]{mathalfa}

\usepackage{multicol}  
\usepackage{epstopdf}
\usepackage{hyperref}

\usepackage{graphicx}
\usepackage{float}

\usepackage{pgfplots}
\pgfplotsset{compat = newest}

\usepackage{multicol}  
\usepackage{epstopdf}

\DeclareUnicodeCharacter{00A0}{~}

\numberwithin{equation}{section}

\newcommand{\N}{\mathbb{N}}
\newcommand{\R}{\mathbb{R}}
\newcommand{\epsi}{\varepsilon}
\renewcommand{\epsilon}{\epsi}

\newcommand{\mi}{\mathcal{i}}
\newcommand{\mj}{\mathcal{j}}
\newcommand{\mij}{\mathcal{ij}}

\renewcommand{\epsilon}{\epsi}

\newcommand{\bK}{\mathbf{K}}

\newcommand{\bA}{\mathbf{A}}

\newcommand{\mM}{\mathcal{M}}
\newcommand{\mN}{\mathcal{N}}
\newcommand{\mNN}{\mathcal{N}\setminus\left\{N\right\}}
\newcommand{\mNZN}{\mathcal{N}\setminus\left\{0,N\right\}}

\DeclareMathOperator{\s}{sign}

\renewcommand{\ss}{\mathrm{s}}

\newtheorem{thm}{Theorem}[section]
\newtheorem{lemma}[thm]{Lemma}
\newtheorem{prop}[thm]{Proposition}

\newtheorem{remark}[thm]{Remark}

\theoremstyle{definition}
\newtheorem{defi}[thm]{Definition}

\title[Opinion formation on evolving network.]{Opinion formation on evolving network. The DPA method applied to a nonlocal cross-diffusion PDE-ODE system}
\author{Simone Fagioli, Gianluca Favre}
\address[Simone Fagioli]{\newline Dipartimento di Ingegneria e Scienze dell’Informazione e Matematica\newline Universit\`a degli Studi dell’Aquila, Via Vetoio 1, 67100 Coppito, L’Aquila, It.}
\email{simone.fagioli@univaq.it}
\address[Gianluca Favre]{\newline Dipartimento di Ingegneria e Scienze dell’Informazione e Matematica\newline Universit\`a degli Studi dell’Aquila, Via Vetoio 1, 67100 Coppito, L’Aquila, It.}
\email{gianluca.favre@univaq.it}
\date{}

\usepackage[final]{showlabels}

\usepackage{marginnote}

\usepackage[backend=bibtex]{biblatex}
\bibliography{references}
\begin{document}

\maketitle
\begin{abstract}
We study a system of nonlocal aggregation cross-diffusion PDEs that describe the evolution of opinion densities on a network. The PDEs are coupled with a system of ODEs that describe the time evolution of the agents on the network. Firstly, we apply the Deterministic Particle Approximation (DPA) method to the aforementioned system in order to prove the existence of solutions under suitable assumptions on the interactions between agents. Later on, we present an explicit model for opinion formation on an evolving network. The opinions evolve based on both the distance between the agents on the network and the 'attitude areas,' which depend on the distance between the agents' opinions. The position of the agents on the network evolves based on the distance between the agents' opinions. The goal is to study radicalization, polarization, and fragmentation of the population while changing its open-mindedness and the radius of interaction.
\end{abstract}

\maketitle

\section{Introduction}
The study of social phenomena through mathematical modelling has gained significant attention in the scientific community, especially in recent decades \cite{BeMaTo,CaFoLo,Ga,MoTa,NaPaTo,PaTobook,St}. The exchange of information on these platforms has sparked research in understanding how social interactions shape the process of opinion formation \cite{AlPaToZa,Be,BoLo,KlShSh,LaMa,SlLa,Sz,YaRoSc}.

In social interactions, the relationships between individuals are often structured as networks that co-evolve with the individuals themselves \cite{ThHaKl}. A prominent example of this is the formation of opinions or norms within social networks, where interactions only occur between connected agents. However, the network connections are dynamic, and this change influences the states of the individuals. For example, opinions can be influenced by connections, such as followers reacting to posts, while individuals tend to follow others with closer opinions.

The network structure of social interactions plays a vital role and is commonly represented using random networks. However, there are two natural levels in examining opinion formation processes on network: the microscopic and macroscopic scales.  The microscopic models have been employed to simulate phenomena such as opinion formation, knowledge networks, social norm formation, and biological transport networks \cite{albi2017continuum,AlPaZa3,benatti2020da,kohne2020role,nigam2018onem,tur2018coevolution}. In considering processes with a huge number of agents, a natural question arise in considering a limit procedure between the two scales. However, the specific details of the network structure can be lost, and only few general characteristics are incorporated into the models \cite{coppini2020law,delattre2016lucon}.

From a mathematical standpoint, it is natural to apply methods from statistical physics or kinetic theory to bridge the gap between microscopic interactions and macroscopic models, see \cite{boud2009,T1}. This approach involves formulating partial differential equations for distributions and using well-established asymptotic methods to simplify the equations and analyse pattern formation \cite{bouchut2000series,cercignani1969mathematical}. These mathematical approaches have successfully explained macroscopic distributions in socio-economic interactions and various aspects of opinion formation and polarization \cite{AlPaToZa,burger2014partial,naldi2010eds,pareschi2014wealth}. Using these tools has significantly contributed to understanding the emergence of macroscopic behaviour from microscopic interactions in a wide range of social phenomena.

\subsection{Modelling motivation: the social context}
During the last two decades, the diffusion of smartphones and the increasing use of social networks have changed how people interact and form their beliefs. There are two main aspects that have been drastically disrupted: the \textbf{number of connections} and the \textbf{frequency of interactions}.

Due to the \textit{hyperconnectivity} of the globalized world, each individual can get in contact with a wide range of opinions. Delocalizing the place of interaction from physical space to the digital realm has destroyed the local cultural bias of interactions. This implies that each individual can come into contact with cultures and ideas they do not know and cannot deeply understand. This aspect has resulted in a change in the epistemic processes, specifically altering the dynamics that govern the formation of beliefs when individuals are exposed to new inputs (news, visual art, songs, posts, tweets, chats, etc.).

Moreover, the increasing amount of inputs and their high frequency create a physical upper bound on the processing capacity of the human brain. An individual walking through a mall cannot process all the inputs coming from screens, speakers, billboards, smartphones, and so on. The same situation occurs while scrolling through social networks or digital social media. As a result, there is a need to filter the inputs, both on the physical side through the network and on the rational side, by selectively processing only a few of them and disregarding others at a peripheral level of thinking.

All these aspects are well understood by scientists from social sciences and social epistemologists (a non-exhaustive list includes references such as \cite{begby_2022, Bernecker_of, Lackey_of, nguyen_2020, Rosenstock_of, O_Connor_of}). However, it is challenging to fit them into a unified mathematical description. In this paper, we propose an approach that mainly focuses on two tools: attitude areas and the Euclidean network. In Section \ref{sec: modelling}, we introduce and simulate a new model to investigate the evolution of the network and opinion distributions for agents interacting on social networks and social media.

\subsection{Mathematical tools}
We consider here a generalisation of the model introduced by Burger in  \cite{Burger_Net,Burger_OF}, where the author derives a kinetic description of the opinion formation process on networks,  studying marginal hierarchies and pass to an infinite limit in several situations. We assume that each agent is described by two variables: the \textit{position on the network} space, and the \textit{opinion density distribution}.
Given the space of opinion-network variables $(\Omega\subset\mathbb{R}) \times \mathbb{R}^d$, the distribution $\rho^\mi(t,x)\in \mathbb{R}^+$ with $(t,x)\in \mathbb{R}^+ \times \Omega$ describes the opinion of the $\mi$-th agent. While, the position on the network is given by $a^\mi(t)\in\mathbb{R}^d$. The model we are going to consider describes the evolution of the opinion and of the agents position in the network. 
Given $\textbf{a}$ and $\bm{\rho}$ the vector of network positions and that of opinion distributions respectively, the evolution is given by the following PDE-ODE system
\begin{subequations}\label{eq:main_intro}
	\begin{align}
		& \partial_t \rho^\mi(t,x) = \partial_x \left( \beta^\mi(\boldsymbol{\rho},\textbf{a};x)  \partial_x \rho^\mi(t,x) \right) -  \partial_x \left( \rho^\mi(t,x) \theta^\mi(\bm{\rho},\textbf{a};x) \right)\,,\\
		&  \partial_t a^\mi(t) = \sum_{\mj\in\mM}\mathbf{V}(\mu_{\rho^\mi}(t), \mu_{\rho^\mj}(t), a^\mij)\,,
	\end{align} 
\end{subequations}
where $\beta^\mi$ and $\theta^\mi$ are defined by
\begin{align}
	\beta^\mi(\bm{\rho},\textbf{a};x) &= \sum_{\mj = 1}^M \int_\Omega \mathbf{A}^\mij(x,y,a^\mij) \rho^\mj(y,t)\,dy\,,\\
	\theta^\mi(\bm{\rho},\textbf{a};x) &=\sum_{\mj = 1}^M \int_\Omega \mathbf{K}^\mij(x,y,a^\mij) \rho^\mj(t,y)\,dy\,,
\end{align}
with $\mathbf{A}^\mij$ and $\mathbf{K}^\mij$ interaction potentials that will be defined later. In \eqref{eq:main_intro} $\mu_{\rho^\mi}$ denotes the first momentum of $\rho^\mi$, i.e.~the mean opinion of the $\mi$-th agent,
%\begin{equation*}
%    \mu_{\rho^\mi}(t) = \bigg(\int_{-\infty}^{+\infty} x \, \rho^\mi(t,x) \, dx \bigg) \bigg( \int_{-\infty}^{+\infty} \rho^\mi(t,x) \, dx \bigg)^{-1}\,.
%\end{equation*}
and we assume that the evolution of agent $\mi$ in the network is influenced by the interaction with agent $\mj$ through the mutual euclidean distance between the agents - defined by $a^\mij(t)$ - and the mean opinions $\mu_{\rho^\mi}$ and $\mu_{\rho^\mj}$. On the other hand, the opinion density distribution obeys a nonlocal aggregation-diffusion equation. Note that both the diffusion mobilities $\beta^\mi$ and the trasport operators $\theta^\mi$ result to be nonlocal operators. Moreover, the dependence of $\beta^\mi$ from the vector $\bm{\rho}$ induces the effect of a cross-diffusion mechanism on the system.

The goal of the paper is twofold. On one hand we are interested in the existence of solution for system \eqref{eq:main_intro} in an appropriate functional setting. On the other hand we want to investigate numerically solutions to system \eqref{eq:main_intro} in order to deduce if it is able to reproduce processes such as polarization, radicalization, fragmentation, and clustering of the population. Similar questions have been investigated  recently by Nugent et al.~\cite{wolfram_of}.

To establish the existence of solution for system \eqref{eq:main_intro}, the authors draw inspiration from the deterministic particle approximation (DPA) developed for similar equations in \cite{FagRad,FaRa,FagTse}. The method and its several modifications date back to the seminal works \cite{GT1,russo1} to demonstrate the convergence of the resulting equation and was applied in several contexts such as traffic flow \cite{diFFRo2017,DiFRo} and local or nonlocal transport equations \cite{DiFFaRa,DiFSt,FagTse}.

The DPA is then used to perform numerical simulations on \eqref{eq:main_intro}. As a numerical scheme it can be connected to moving mesh schemes applied in various contexts, such as diffusion problems, computational fluid dynamics, and scalar conservation laws \cite{BHJ11,BHR96,BHR09,CHR02,CHR03,SMR01}. However, it should be noted that these schemes also have limitations in one-dimensional applications, similar to scalar conservation laws.

The paper is structured as follows. In Section \ref{sec: 2} we first introduce the rigorous Deterministic Particle Approximation (DPA) we are dealing with. We introduce some preliminaries on optimal transportation theory together with the main assumptions. We conclude the section with the statement of the main result in Theorem \ref{th: main thm}. Section \ref{sec: 3} is devoted to the proof of the main theorem by providing fundamental a priori estimates that allow to deduce convergence of proper reconstructed piecewise constant densities to weak solutions. Finally, in Section \ref{sec: modelling} we make use of DPA numerical scheme in order to simulate an explicit model for opinion formation on evolving network. The aim is about studying radicalization, polarization, and fragmentation of the population while changing its open mindedness and the radius of interaction.

\section{Rigorous formulation, assumptions and main result}\label{sec: 2}
\subsection{Deterministic Particle Approximation (DPA)} We begin this section with the rigorous formulation of the particle evolution already sketched in the Introduction. We consider in $\mathbb{R}^d$ a network of $M$ nodes and we locate an agent $a^\mi\in\mathbb{R}^d$ with $\mathcal{i}\in \mM=\{1,...,M\}$ in each node. Assume that each agent may have opinion ranging on a finite set $\Omega \subset \mathbb{R}$, without loss of generality we consider $\Omega = [-1,1]$. To each agent we associate a finite \emph{opinion strength} $\sigma^\mi$ and an initial opinion density $\bar{\rho}^\mi(x)\in L^1(\Omega)$ such that
\begin{equation*}
	\sigma^\mi = \int_\Omega \bar{\rho}^\mi (y)\,dy\,, \qquad \forall\mathcal{i}\in \mM\,.
\end{equation*}
Given $N\in\mathbb{N}$, we consider the strength fractions $\sigma^\mi_N = \sigma^\mi / N$, and we introduce for each $\mi \in \mM$ the $\{\bar{x}_k^\mi\}$ partition of $\Omega$ with $k\in\mN=\{0,...,N\}$ given by
\begin{equation}\label{part_ini}
	\begin{aligned}
		&\bar{x}_0^\mi = -1 \,,\\
		&\bar{x}_k^\mi = \inf\left\{ x \in \Omega \, : \, \int_{x_{k-1}}^x \bar{\rho}^\mi (y)\,dy = \sigma^\mi_N \right\} \,, \qquad k\in \mNZN\\
		&\bar{x}_N^\mi = 1 \,.
	\end{aligned} 
\end{equation}
Note that $\bar{x}_k^\mi<\bar{x}_{k+1}^\mi$, for any $\mi\in\mM$ and $k\in\mNN$. This procedure allow to associate to each agent a finite number of time-evolving opinions $x^\mi_k(t)$. Assume that initially all the nodes $\bar{a}^\mi = a^\mi(t=0)$ are located in a certain smooth and bounded domain $\Lambda\in\mathbb{R}^d$. We then let the nodes evolve in time depending on the distances $a^\mij$ between the agents $a^\mi$ and $a^\mj$, and the mean opinion of the agents. 

We define the \emph{discrete opinion densities} for the $\mi$-th agent as
\begin{equation}\label{eq:disc_den}
	\rho^\mi_k(t) = \frac{\sigma^\mi_N}{|I_k^\mi(t)|},\quad\mbox{ with }\quad I_k^\mi(t)=[x^\mi_k(t),x^\mi_{k+1}(t))
\end{equation}
with $k\in\mNN$, and the \textit{discrete mean opinions} by
\begin{equation}
	\label{def: average x}
	\mu_{x^\mi}^{\mN}(t) = \frac{1}{N+1} \sum_{k\in\mN} x^\mi_k(t)\,,
\end{equation}
where $a^\mij_\mN$ is the euclidean distance between the agents after the opinion discretization.
In the following we may denote with $\textbf{x}^{\mi,\mN}(t):=\left(x_0^\mi(t),\ldots,x_N^\mi(t)\right)$, for all $\mi\in\mM$.

Thus we consider the following system of ODEs
\begin{subequations}\label{eq:DPA}
	\begin{align}
		\label{eq: evol x}
		&\dot{x}^\mi_k(t) = \frac{\beta^\mi_k}{\sigma_N^\mi} \left(\Phi^\mi(\rho^\mi_{k-1}) - \Phi^\mi(\rho^\mi_{k}) \right) + \theta^\mi_k \\
		\label{eq: evol agent}    
		&\dot{a}_\mN^\mi (t) = \sum_{\mj\in\mM}\mathbf{V}(\mu_{x^\mi}^{\mN}, \mu_{x^\mj}^{\mN}, a_\mN^{\mij}), 
	\end{align}
\end{subequations}
for $k\in\mNZN$ and $\mi \in \mM$, endowed with the \emph{boundary conditions}
\begin{equation}
	\label{eq: bound}    
	\dot{x}^\mi_0(t) = 0\quad, \quad
	\dot{x}^\mi_{N}(t) = 0 \qquad \text{for all }\mi \in \mM,
\end{equation}
and initial conditions
\begin{equation}
	\label{eq: initial}    
	x^\mi_k(0) = \bar{x}^\mi_k\quad , \quad
	a_\mN^\mi(0) = \bar{a}^\mi \qquad \text{for all }\mi \in \mM,\,k\in\mN.
\end{equation}

In \eqref{eq:DPA} have denoted with $\beta^\mi_k(t)$ the \emph{discrete diffusion mobilities}
\begin{equation}
	\label{eq: def beta}
	\beta^\mi_k(t) =  \sum_{\mj\in\mM} \sum_{l\in\mN} \sigma^\mj_N\mathbf{A}^{\mij}(x^\mi_k,x^\mj_l,a_\mN^\mij)\,,
\end{equation} 
and with $\Phi^\mi$ the nonlinear diffusion for the agent $\mi$. Functions $\theta^\mi_k(t)$ describe the \emph{discrete transports} and are given by
\begin{equation}
	\label{eq: def theta}    
	\theta^\mi_k(t) =  \sum_{\mj\in\mM} \sum_{l\in\mN} \sigma^\mj_N\mathbf{K}^{\mij}(x^\mi_k,x^\mj_l,a_\mN^\mij)\,.
\end{equation}

\subsection{Preliminaries and assumptions} We now present some tools from optimal transport that will be useful in the sequel. The Wasserstein distance is the right notion of distance for the opinions since it allows to measure distances between measures (densities) with same mass.  For a fixed mass $\sigma>0$, we consider the space
\begin{equation*}
	\mathfrak{M}_\sigma = \bigl\{\mu \hbox{ Radon measure on $\R$ } \colon \mu\ge 0 \text{ and }\mu(\R)=\sigma \bigr\}.
\end{equation*}
Given $\mu\in \mathfrak{M}_\sigma$, we introduce the pseudo-inverse function $X_\mu \in L^1([0,\sigma];\R)$ as
\begin{equation}\label{eq:pseudoinverse}
	X_\mu(z) = \inf \bigl\{ x \in \R \colon \mu((-\infty,x]) > z \bigr\}.
\end{equation}
In particular, if $\sigma=1$, then $\mathfrak{M}_1$ is the set of non-negative probability densities on $\R$ and it is possible to consider the one-dimensional \emph{$1$-Wasserstein distance} between each pair of densities $\rho_1,\rho_2\in \mathfrak{M}_1$. As shown in \cite{CT}, in the one dimensional setting the \emph{$p$-Wasserstein distance} can be equivalently defined in terms of the $L^1$-distance between the respective pseudo-inverse mappings as
\[
d_{W^p}(\rho_1,\rho_2) = \|X_{\rho_1}-X_{\rho_2}\|_{L^p([0,1];\R)}.
\]
For generic $\sigma >0$, we recall the definition for the \emph{scaled $1$-Wasserstein distance} between $\rho_1,\rho_2\in \mathfrak{M}_\sigma$ as
\begin{equation}\label{eq:wass_equiv0}
	d_{W^1_\sigma}(\rho_1,\rho_2)= \|X_{\rho_1}-X_{\rho_2}\|_{L^1([0,\sigma];\R)},
\end{equation}
We refer to \cite{AGS,S,V1} for a complete presentation of the subject.

We assume that the \textbf{initial densities} are under the following assumptions:
\begin{itemize}
	\item[(In1)] $\bar{\rho}^\mi \in BV(\Omega; \R^+)$ with $\|\bar{\rho}^\mi\|_{L^1(\Omega)} = \sigma^\mi$, for some $\sigma^\mi>0$,
	\item[(In2)] there exist $m^\mi,M^\mi >0$ such that $m^\mi \leq \bar{\rho}^\mi(x) \leq M^\mi$ for every $x \in \Omega$.
\end{itemize}

We now introduce the assumptions for the \textbf{diffusive} and \textbf{transport operators}.
\begin{itemize}
	\item[($\bA$)] We assume that $\bA^\mij : \Omega\times\Omega\times\R\to \R$ is a $C^1$ non-negative function w.r.t.~the first variable and for all pairs $(\mi,\mj)$ it exists $c_\bA > 0$ such that
	\begin{equation*}
		|\bA^\mij(x^i_{k^*}, x^j_{s^*}, a^\mij) - \bA^\mij(x^i_k, x^j_s,a^\mij)| \le c_\bA \left( |x^i_{k^*} - x^i_k| + |x^j_{s^*} - x^j_s| \right)\,,
	\end{equation*}
	and it exists $c_{1,\bA} > 0$ such that
	\begin{equation*}
		|\partial_1\bA^\mij(x^i_{k^*}, x^j_{s^*}, a^\mij)-\partial_1\bA^\mij(x^i_k, x^j_s,a^\mij)|\leq c_{1,\bA}|x^j_{s^*} - x^j_s| ,
	\end{equation*}
	for all $(k,k^*)$ and $(s,s^*)$ pairs of indexes in $\mN\times\mN$, where $\partial_1\bA^\mij$ denotes the derivatives with respect to the first entrance.
	\item[($\bK$)] We assume that $\bK^\mij : \Omega\times\Omega\times\R\to \R$ is bounded, continuous, and for all pairs $(\mi,\mj)$ it exists $c_\bK > 0$ such that
	\begin{equation*}
		|\bK^\mij(x^i_{k^*}, x^j_{s^*},a^\mij) - \bK^\mij(x^i_k, x^j_s,a^\mij)| \le c_\bK \left( |x^i_{k^*} - x^i_k| + |x^j_{s^*} - x^j_s|\right)\,,
	\end{equation*}
	for all $(k,k^*)$ and $(s,s^*)$ pairs of indexes in $\mN\times\mN$. We further assume that $$\mathbf{K^{\mi\mi}}(x,x,a^{\mi\mi})=0$$.
	\item[(Dif)] $\Phi^\mi: [0,\infty) \to \R$ is a nondecreasing Lipschitz function, with $\Phi^\mi(0)=0$.
	\item[($\mathbf{V}$)] The network velocity $\mathbf{V}$ is a $C^1$ bounded function on $\Omega \times \Omega \times \R^{+}$.
\end{itemize}

\subsection{Continuous reconstruction and main result}
Given the preliminaries assumptions, we give the definition of weak solutions to equation \eqref{eq:main_intro} together with the statement of the main result. 

By setting $\Omega_T = [0,T]\times\Omega$ and $\partial\Omega_T=[0,T]\times\left\{-1,1\right\}$, and considering $\bar{\rho}^\mi\in L^1\cap L^\infty(\Omega)$ and $\bar{a}^\mi\in\Lambda\subset\R^d$, we are going to deal with the following PDE-ODE system
\begin{equation} \label{eq:IBVP}
	\begin{dcases}
		\partial_t \rho^\mi = \partial_x \left( \beta^\mi(\bm{\rho},\textbf{a};x)  \partial_x \Phi^\mi(\rho^\mi )\right) - \partial_x \left( \rho^\mi\theta^\mi(\bm{\rho},\textbf{a};x)\right),& (t,x)\in \Omega_T,\\
		\partial_t a^\mi(t) = \sum_{\mj\in\mM}\mathbf{V}(\mu_{\rho^\mi}(t), \mu_{\rho^\mj}(t), a^\mij),& t\in [0,T],\\
		\beta^\mi(\bm{\rho},\textbf{a};x)  \partial_x \Phi^\mi(\rho^\mi(t,x) )-\rho^\mi(t,x)\theta^\mi(\bm{\rho},\textbf{a};x)=0,& (t,x)\in \partial\Omega_T,\\
		\rho^\mi(0,\cdot) = \bar{\rho}^\mi, & x\in\Omega,\\
		a^\mi(0)=\bar{a}^\mi,
	\end{dcases}\, 
\end{equation}
for all $\mi\in\mM$, where the bold notation refers to the vectors $\bm{\rho} = (\rho^1,...,\rho^M)$ and the set $\textbf{a}=(a^1,...,a^M)$. 

We state the notion of weak solution for the system \eqref{eq:IBVP} as follows 
\begin{defi}[Weak solution]\label{def:weak}
	We say that the couple $\left(\bm{\rho},\textbf{a}\right)$ is weak solution of \eqref{eq:main_intro} in the formulation \eqref{eq:IBVP} if 
	\begin{itemize}
		\item $\rho^\mi \in L^\infty  \cap BV(\Omega_T)$, with $\rho^\mi(0,\cdot)= \bar{\rho}^\mi$, for all $\mi\in \mM$
		\item $a^{\mi}\in C^2\left(\left[0,T\right];\R^d\right)$,
	\end{itemize}
	and taken $\zeta\in C_0^{\infty}(\Omega_T)$, for all $\mi\in \mM$ it satisfies
	\begin{equation}\label{weak_rho}
		\begin{aligned}
			\int_{\Omega_T} \rho^\mi(t,x) &\partial_t\zeta(t,x)+\rho^\mi(t,x)\theta^\mi(\bm{\rho},\textbf{a};x) \partial_x \zeta(x)\, dx\,dt \\
			+\int_{\Omega_T}&\Phi^\mi\left(\rho^\mi(t,x)\right) \left( \partial_x\beta^\mi(\bm{\rho},\textbf{a};x)   \partial_x \zeta(t,x)+\beta^\mi(\bm{\rho},\textbf{a};x)   \partial_{xx} \zeta(t,x)\right)\, dx=0\,\end{aligned}
	\end{equation}
	and
	\begin{equation}\label{mild_a}
		a^{\mi}(t) = \bar{a}^{\mi}+\sum_{\mj\in\mM
		}\int_0^t  \mathbf{V}\left(\mu_{\rho^\mi}(\tau), \mu_{\rho^\mj}(\tau), a^\mij(\tau)\right)\,d\tau
	\end{equation}
	for all $t\in [0,T]$.
\end{defi}

Given the discrete opinion densities defined in \eqref{eq:disc_den} we consider the following piece wise constant density reconstructions
\begin{equation}\label{eq:piec_discden}
	\rho^{\mi,\mN}(t,x) = \sum_{k\in\mNN}\rho^\mi_k(t)\, \chi_{I_k^\mi(t)}(x)\,.
\end{equation}

The main result of the paper reads as follows
\begin{thm}\label{th: main thm}
	Given $M\in\N$ and $T>0$ fixed, consider $\bA^{\mij}$, $\bK^{\mij}$, $\mathbf{V}$, and $\Phi^\mi$ under assumptions $(\bA)$, $(\bK)$, $(\textbf{V})$, and $(Dif)$ respectively, for all $\mi,\mj\in \mM$. Let $\bar{\rho}^\mi:\Omega \to \R$ under assumptions $(In1)$ and $(In2)$ and $\bar{a}^\mi\in\Lambda\subset\R^d$ for all $\mi\in\mM$. Then, for all $\mi\in \mM$, when $N\to \infty$ the density $\rho^{\mi,\mN}$ introduced in equation \eqref{eq:piec_discden} converges (up to subsequence) strongly to $\rho^\mi\in L^\infty\cap BV (\Omega_T)$ and the solution to equation \eqref{eq: evol agent} $a_\mN^{\mi}$ converges to $a^{\mi}\in C^2\left(\left[0,T\right];\R^d\right)$  where the couple $\left(\bm{\rho},\textbf{a}\right)$ is a solution to equation \eqref{eq:main_intro} in the sense of Definition \ref{def:weak}. 
\end{thm}

To avoid lack of notation, we highlight that now on while considering the limit $N\to\infty$ we refer to the limit to infinity of the cardinality of the set of indexes $\mN=\{0,...,N\}$. 

\section{Proof of the main result}\label{sec: 3}
\subsection{Basic estimates} We start providing some fundamental estimates that allow to deduce the well-posedness of \eqref{eq:DPA} and of the discrete densities \eqref{eq:disc_den}. Let us recall that in \eqref{eq:disc_den} we introduced the intervals
\begin{equation}\label{qe:interval}
	I_k^\mi(t) = \left[{x}^\mi_{k}(t),{x}^\mi_{k+1}(t)\right),\quad |I_k^\mi|(t)=|{x}^\mi_{k+1} - {x}^\mi_{k}|,
\end{equation}
for $\mi\in\mM$ and $k\in\mNN$. The first step is to prove that such intervals are well-defined. 
We start with the following auxiliary lemma, that follows directly from Assumption $(\bK)$
\begin{lemma}\label{lemma: kappa sum}
	With the setting of the Main Theorem \ref{th: main thm}, and referring to the previous definitions, given $\mathbf{K^\mij}$ under Assumption $(\bK)$, then it exists $C>0$ such that the following inequality holds
	\begin{equation}\label{boundtheta}
		|\theta_{k+1}^i(t) - \theta_k^i(t)| \le C |x^i_{k+1}(t) - x^i_{k}(t)| \quad \forall \mi \in \mM, \, \forall t \in [0,T]\,,
	\end{equation}
	with $C= c_\bK \sigma^\mM $, where $\sigma^\mM = \sum_{j\in\mM} \sigma^\mj$.
\end{lemma}

\begin{lemma}[Ordering preservation]\label{lem: bound from below}
	Assume $\bA^{\mij}$ and $\bK^{\mij}$ under assumptions $(\bA)$ and $(\bK)$ respectively for all $\mi,\mj\in \mM$.
	Let consider the DPA system described by \eqref{eq:DPA} with initial conditions $\bar{x}^\mi$ constructed in \eqref{part_ini}, for $\mi\in \mM$, and a finite time $T>0$. Then, for all $t\in [0,T)$ there is a positive constant $\mu$ independent from $\mN$ - and so from $N$ too - such that the distance between two adjacent opinions $x^\mi_k, x^\mi_{k+1} \in C[0,T]$ is bounded from below by 
	\begin{equation*}
		\left(x^\mi_{k+1} - x^\mi_k\right)(t) \ge \underset{\mi,k}{\min}(\bar{x}^\mi_{k+1} - \bar{x}^\mi_k) \, e^{-\mu T}, 
	\end{equation*}
	for all  $k\in\mNN$ and $\mi \in \mM$.
\end{lemma}

\begin{proof}
	Given $T>0$ and $\mi\in\mM$, we define $\tau_1$ as 
	\begin{equation*}
		\tau_1 = \inf \big\{ s \in (0,T) \,:\, \exists\, k \in \mN\backslash\{N\} \text{ s.t. } (x^{\mi}_{k+1} - x^{\mi}_{k} )(s) = (\bar{x}^{\mi}_{k+1} - \bar{x}^{\mi}_{k})\, e^{-\mu s}    \big\}\,,
	\end{equation*}
	then the same index $k$ corresponds also too the one of the maximum $\rho^\mi_k$ at time $\tau_1$ because $I^i_k$ is the minimum interval of the $\mN$ partition of $\Omega$ for the $i$-th agent at time $\tau_1$.
	
	At this point, let assume that exists $\tau_2 \in (\tau_1 , T)$ such that $$\left(x^\mi_{k+1} - x^\mi_k\right)(s) < (\bar{x}^\mi_{k+1} - \bar{x}^\mi_k) e^{-\mu s} \quad \forall s \in(\tau_1,\tau_2)\,.$$ We show that the existence of $\tau_2$ would bring to a contradiction. Let consider the evolution of the interval $I^i_k$
	\begin{align*}
		\frac{d}{dt} \left[ e^{\mu t} \left(x^\mi_{k+1} - x^\mi_k\right)(t)\right]_{|t=\tau_1} &= e^{\mu \tau_1 } \frac{\beta^\mi_{k+1}(\tau_1) }{\sigma_N^\mi}\left[ \Phi^\mi(\rho^\mi_{k}(\tau_1) ) - \Phi^\mi(\rho^\mi_{k+1}(\tau_1) ) \right] + e^{\mu \tau_1} \theta^\mi_{k+1}(\tau_1)  \\
		&\quad - e^{\mu \tau_1} \frac{\beta^\mi_{k}(\tau_1) }{\sigma_N^\mi}\left[ \Phi^\mi(\rho^\mi_{k-1}(\tau_1) ) - \Phi^\mi(\rho^\mi_{k}(\tau_1) )\right] - e^{\mu \tau_1} \theta^\mi_{k}(\tau_1)  \\
		&\quad + \mu e^{\mu \tau_1} (x^\mi_{k+1}(\tau_1)  - x^\mi_k(\tau_1) )\,.
	\end{align*}
	At time $t=\tau_1$, as highlighted before, by construction of the discrete densities in \eqref{eq:disc_den} we have $\rho^\mi_k \ge \rho^\mi_l$ for all $l\ne k$. 
	
	Then, the monotonicity of $\Phi^\mi$ gives
	\begin{align*}
		\frac{d}{dt} \left[ e^{\mu t} \left(x^i_{k+1} - x^i_k\right)(t)\right]_{|t=\tau_1} &\ge  e^{\mu \tau_1}\left[ \mu\, \left(x^\mi_{k+1}(\tau_1)  - x^\mi_k(\tau_1) \right)+\left(\theta^\mi_{k+1}(\tau_1) - \theta^\mi_{k}(\tau_1)\right) \right].
	\end{align*}
	Thanks to \eqref{boundtheta} we get 
	\begin{align*}
		\frac{d}{dt} \left[ e^{\mu t} \left(x^\mi_{k+1} - x^\mi_k\right)(t)\right]_{|t=\tau_1} &\ge e^{\mu \tau_1}\, (\mu - C) (x^\mi_{k+1}(\tau_1) - x^\mi_{k}(\tau_1)) \ge 0\,,  
	\end{align*}
	which holds while choosing $\mu \ge C$.
	At this point we fix $t^*\in (\tau_1,\tau_1+\delta < \tau_2)$ with $\delta$ as small as we want, and we get
	\begin{equation*}
		e^{\mu t^*} \left(x^\mi_{k+1} - x^\mi_k\right)(t^*) = e^{\mu \tau_1} \left(x^\mi_{k+1} - x^\mi_k\right)(\tau_1) + \int_{\tau_1}^{t^*} \frac{d}{ds} \left[ e^{\mu t} \left(x^\mi_{k+1} - x^\mi_k\right)(s)\right] \, ds \,,
	\end{equation*}
	due to the positiveness of the last term we get the wished result which show the absurd,
	\begin{equation*}
		\left(x^\mi_{k+1} - x^\mi_k\right)(t^*) \ge e^{- \mu t^*} \left(\bar{x}^\mi_{k+1} - \bar{x}^\mi_k\right)\,,
	\end{equation*}
	this proves that $I^i_k(t)$ cannot decrease faster than $I^i_k(0) e^{- \mu t}$. Nevertheless, this does not deny the existence of an index $\tilde{k}$ such that the interval $I^i_{\tilde{k}}(t)$, satisfying $I^i_{\tilde{k}} (\tau_1) > I^i_{\tilde{k}}(0)$, decreases faster than $I^i_k(t)$. This means that could exists $\tilde{\tau}_1 > \tau_1$ for which $I^i_{\tilde{k}}(\tilde{\tau}_1) = I^i_{\tilde{k}}(0)$, with $\tilde{\tau}_1 < T$. At this point we should prove that there exists $\tilde{\mu}$ such that $I^i_{\tilde{k}}(t) \ge e^{-\tilde{\mu}t} I^i_{\tilde{k}}(0)$ for all $t \in (\tilde{\tau}_1 , T)$. To prove it, we repeat the same procedure explained before but defining $\tilde{\tau}_1$ considering the set of indexes $\mN / \{k\}$. The final exponential rate will be the largest $\mu$ among those considered.
	
	%\subsubsection*{Caveat}
	Let also consider the case with $k$ not unique, i.e.~for $\tau_1$ there are several intervals satisfying the definition of $\tau_1$, the set of these indexes is denoted by $\mathcal{J}=\{ k_j \}$. If there is at least one $k_{j^*}$ not adjacent to other indexes of $\{k_j\}$, then we take that index and we follow the proof above. In the event that there are three indexes of $\{k_j\}$ in a row, we take without distinction that one with the fastest decrease in time of the interval $I^i_{k_j}$. At this point we are back to the steps showed above. This concludes the proof.

\end{proof}

\begin{remark}
	A similar procedure of the one in Lemma \ref{lem: bound from below} allows to produce the upper bound on discrete opinions
	\begin{equation}\label{eq:upbound opinion}
		\left(x^\mi_{k+1} - x^\mi_k\right)(t) \le \underset{k}{\max}(\bar{x}^\mi_{k+1} - \bar{x}^\mi_k) \, e^{CT} \qquad \forall t\in[0,T]\,.
	\end{equation}
	This estimate, together with the one in Lemma \ref{lem: bound from below}, allows to deduce the following bounds on the discrete densities in \eqref{eq:disc_den}
	\begin{equation}\label{eq:bound or rhoik}
		m^\mi e^{-CT}\leq \rho_{k}^\mi\leq M^\mi e^{\mu T}\,, \mbox{ for all }\,\mi\in\mM,\,\,k\in\mNN,
	\end{equation}
	with $m^\mi$ and $M^\mi$ in assumption (In2).
\end{remark}

\begin{lemma}[Velocity boundedness]\label{lem:bound_velocity}
	Assume $\bA^{\mij}$ and $\bK^{\mij}$ respectively under assumptions $(\bA)$ and $(\bK)$ for all $(\mi,\mj)\in \mM\times\mM$.
	Given $\mi\in \mM$, and a finite time $T>0$, let consider the DPA system described by equation\eqref{eq:DPA} with initial conditions $\bar{x}^\mi$ constructed as in equation \eqref{part_ini}. Then,
	\begin{equation*}
		\sup_{t\in\left[0,T\right]}\|\dot{\mathbf{x}}^{\mi,\mN}(t)\|_{\infty}<+\infty, \,\mbox{ for all }\, \mi\in\mM.
	\end{equation*}
\end{lemma}
\begin{proof}
	Using the equation for the evolution of the partitioning we get
	\begin{align*}
		\frac{1}{2}\frac{d}{dt} \left(\dot{x}^\mi_k(t)\right) ^2 &\le  x_k^\mi \left( \frac{\beta^\mi_{k}}{\sigma_N^\mi} \left(\Phi^\mi(\rho^\mi_{k-1}) - \Phi^\mi(\rho^\mi_{k}) \right) -\theta^\mi_k\right)\, \\
		&\le |x_k^\mi|\left( \frac{1}{\sigma_N^\mi}\|\beta^\mi_{k}\|_{\infty} Lip\left(\Phi^\mi\right)\Big|  \rho^\mi_{k} - \rho^\mi_{k-1}\Big| + \Big|\theta^\mi_k\Big|\right)\\
		\leq&c_1|x_k^\mi|^2+c_2\Big|  \rho^\mi_{k} - \rho^\mi_{k-1}\Big|^2+c_3\Big|\theta^\mi_k\Big|^2,
	\end{align*}
	for some constants $c_1,c_2,c_3>0$. Then, the thesis follows from equation \eqref{eq:bound or rhoik}, assumption $(\bK)$, and Gronwall's Lemma.
\end{proof}

We now have all the tools needed to prove the convergence in some strong sense for the piece wise constant densities of equation \eqref{eq:piec_discden}. We based our strategy on the ones proposed in the context of DPA, see for instance the proofs in \cite{FagRad,FagTse}, that are using the generalised Aubin-Lions lemma version in \cite{RoSa}, that we report here in a simplified version adapted to our setting.

\begin{thm}\label{Aubin Lions}
	
	Let $T >0$ be fixed, and $\rho^\mN(t,\cdot): [a,b] \to \mathbb{R}$ be a sequence of non negative probability densities for every $t \in [0,T]$ and for every $N \in \mathbb{N}$, where $\mN = \{0,...,N\}$. Moreover, assume that $\| \rho^\mN(t,\cdot )\|_{L^\infty} \leq M $ for some constant $M$ independent on $t$ and $N$. If 
	\begin{itemize}
		\item[I)] $\sup_N \int_0^T TV[\rho^\mN(t,\cdot)]dt < \infty$,
		\item[II)] $d_{W^1}(\rho^{\mN}(t,\cdot),\rho^{\mN}(s,\cdot)) < C|t-s|$ for all $t,s \in [0,T]$, where $C$ is a positive constant independent on $N$,
	\end{itemize}
	then $\rho^\mN$ is strongly relatively compact in $L^1([0,T]\times[a,b])$.
\end{thm}  
The result reads as follows
\begin{prop}\label{comp}
	Let $\rho^{\mi,\mN}$ be defined as in equation\eqref{eq:piec_discden} for $\mi\in\mM$. Then there exist $\rho^\mi \in L^1 \cap L^{\infty}(\Omega_T)$ such that $\|\rho^{\mi,\mN} - \rho^\mi\|_{L^1} \to 0$ as $N \to \infty$.
\end{prop}
\begin{proof}
	The proof reduces to the application of Theorem \ref{Aubin Lions}, in order to show that we can apply that result we first prove that
	\begin{equation}\label{BVest}
		\sup_N \int_0^T TV[\rho^{\mi,\mN}(t,\cdot)]\,dt \le \infty\,.
	\end{equation}
	To show this result we look for a Gr\"onwall type inequality. The first step consists of proving the Lipschitz continuity in time of $t \to TV[\rho_t^{\mi,\mN}(t,\cdot)]$. From the boundedness of $\frac{\beta^\mi_{k}}{\sigma_N^\mi}$, and by Lemma \ref{lem: bound from below}, it follows
	\begin{align*}
		TV[\rho^{\mi,\mN}(t,\cdot)] &= \rho_0^\mi(t) + \sum_{k\in\mNZN} |\rho_k^\mi(t) - \rho_{k-1}^\mi(t)| + \rho_{N-1}^\mi(t) \\
		&\le \frac{h\, e^{\mu T}}{\underset{k,\mi}{\min} |\bar{x}^\mi_{k+1} - \bar{x}^\mi_{k}|^2 } \max_{k,\mi}\frac{\beta^\mi_{k}}{\sigma_N^\mi} |t-s|
	\end{align*}
	for all $s$ and $t$ satisfying $0\le s\le t \le T$. We now consider the time derivative
	\begin{equation*} 
		\frac{d}{dt}TV[\rho^{\mi,\mN}(t,\cdot)] = \dot{\rho}_0^\mi(t) + \sum_{k\in\mNZN}\s\big({\rho}_k^\mi(t) - {\rho}_{k-1}^\mi(t)\big) \, \left(\dot{\rho}_k^\mi(t) - \dot{\rho}_{k-1}^\mi(t)\right) + \dot{\rho}_{N-1}^\mi(t) \,.
	\end{equation*}
	Rearranging the sum and defining the operator 
	\begin{equation*}
		\ss_{k} = \begin{cases}
			1- \s\big({\rho}_1^\mi(t) - {\rho}_{0}^\mi(t)\big) & k=0, \\
			\s\big({\rho}_k^\mi(t) - {\rho}_{k-1}^\mi(t)\big) - \s\big({\rho}_{k+1}^\mi(t) - {\rho}_{k}^\mi(t)\big),&k=1,\ldots,N-2,\\
			1+\s\big({\rho}_{N-1}^\mi(t) - {\rho}_{N-2}^\mi(t)\big), & k=N-1,
		\end{cases}
	\end{equation*}
	with we can rewrite the previous equation as
	\begin{equation}
		\label{eq: 2}
		\frac{d}{dt}TV[\rho_t^{\mi,\mN}(t,\cdot)] = \ss_{0}\dot{\rho}_0^\mi(t) + \sum_{k=1}^{N-2} \ss_{k} \dot{\rho}_k^\mi(t) + \ss_{N-1}\dot{\rho}_{N-1}^\mi(t) \,.
	\end{equation}
	
	We now compute
	
	\begin{align*}
		\dot{\rho}_{k}^\mi & = -\frac{\rho_{k}^\mi}{|I_k^\mi|} \Big( \underbrace{\frac{\beta^\mi_{k+1}}{\sigma_N^\mi} (\Phi^\mi\left(\rho_k^\mi\right) - \Phi^\mi\left(\rho_{k+1}^\mi\right)) - \frac{\beta^\mi_{k}}{\sigma_N^\mi} (\Phi^\mi\left(\rho_{k-1}^\mi\right) - \Phi^\mi\left(\rho_{k}^\mi\right))}_{:=B_k^\mi} + \underbrace{(\theta_{k+1}^\mi - \theta_{k}^\mi)}_{:=|I_k^\mi|\Theta_k^\mi} \Big) \\&= -\frac{{\rho}_{k}^\mi}{|I_k^\mi|}B_k^\mi+ \rho_{k}^\mi \Theta_k^\mi.
	\end{align*}
	
	At this point, we show that the terms involving the diffusion are always negative, i.e.
	\[
	-\ss_k\frac{{\rho}_{k}^\mi}{|I_k^\mi|}B_k^\mi\leq 0\quad k=2,\ldots,N-2.
	\]
	In order to prove this statement we should distinguish different cases. First, we observe that $\ss_k$ is always zero if $\rho_{k-1}^\mi<\rho_{k}^\mi<\rho_{k+1}^\mi$ or $\rho_{k+1}^\mi<\rho_{k}^\mi<\rho_{k-1}^\mi$. In the other two cases, namely $\rho_{k}^\mi<\rho_{k-1}^\mi$ and $\rho_{k}^\mi<\rho_{k+1}^\mi$, or $\rho_{k-1}^\mi,\,\rho_{k+1}^\mi<\rho_{k}^\mi$ and $\rho_{k-1}^\mi<\rho_{k}^\mi$, the monotonicity of $\Phi^\mi$ implies respectively $\ss_k\leq 0$ and $B_k^\mi\leq 0$, and $\ss_k\geq 0$ and $B_k^\mi\geq 0$, hence the negativity of the diffusion contribution is proved.
	
	Concerning the term involving $\Theta_k^\mi$ we can rearrange the sum as follows
	\begin{align*}
		\sum_{k=1}^{N-2}\ss_k\rho_k^\mi\Theta_k^\mi  = & \s(\rho_1^\mi-\rho_0^\mi)\rho_1^\mi\Theta_1^\mi-\s(\rho_{N-1}^\mi-\rho_{N-2}^\mi)\rho_{N-2}^\mi\Theta_{N-2}^\mi\\
		&+ \sum_{k=2}^{N-2}\s(\rho_{k}^\mi-\rho_{k-1}^\mi)(\rho_{k}^\mi-\rho_{k-1}^\mi)\Theta_{k}^\mi\\
		&+ \sum_{k=2}^{N-2}\s(\rho_{k}^\mi-\rho_{k-1}^\mi)(\Theta_{k}^\mi-\Theta_{k-1}^\mi)\rho_{k-1}^\mi.
	\end{align*}
	Observing that $|\Theta_k^\mi|\leq C$ because of equation \eqref{boundtheta}, and that
	\begin{align*}
		\left|\theta_{k}^\mi-\theta_{k-1}^\mi\right|\leq &\left|\frac{1}{|I_k^\mi|}\left[\left(\theta_{k+1}^\mi-\theta_{k}^\mi\right)-\left(\theta_{k}^\mi-\theta_{k-1}^\mi\right)\right]\right|\\&+\left|\left(\frac{1}{|I_k^\mi|}-\frac{1}{|I_{k-1}^\mi|}\right)\left(\theta_{k}^\mi-\theta_{k-1}^\mi\right)\right|\\
		\leq & C\frac{\rho_k^\mi}{\sigma_N^\mi}\left(|I_k^\mi|^2+|I_{k-1}^\mi|^2+\left||I_{k}^\mi|-|I_{k-1}^\mi|\right|\right)+\frac{C}{\sigma_N^\mi}\left|\rho_{k}^\mi-\rho_{k-1}^\mi\right||I_{k-1}^\mi|,
	\end{align*}
	we can bound
	\begin{align*}
		\left|\sum_{k=1}^{N-2}\ss_k\rho_k^\mi\theta_k^\mi\right|  \leq & 2M^\mi C\left(1+|\Omega|\right)+ 3C\, TV[\rho^{\mi,\mN}(t)].
	\end{align*}
	We can finally estimate
	\begin{align*}
		\frac{d}{dt}TV[\rho_t^{\mi,\mN}(t,\cdot)] \leq & 2C\frac{M^\mi}{m^\mi}\left(\rho_0^\mi(t) +\rho_{N-1}^\mi(t)\right)+ 2M^\mi C\left(1+|\Omega|\right)+ 3C\, TV[\rho^{\mi,\mN}(t)] \,,
	\end{align*}
	and thus equation \eqref{BVest} follows by Gr\"onwall type argument.
	
	We now prove that the second requirement of Theorem \ref{Aubin Lions} holds, namely there exists a positive constant $C$ such that
	\begin{equation}\label{w_1bound}
		d_{W^1}\big(\rho^{\mi,\mN}(t,\cdot),\rho^{\mi,\mN}(s,\cdot)\big) \le C |t-s| \qquad \forall\,  s,t\in[0,T]\,.
	\end{equation}
	In order to do this, we use the isometry of equation \eqref{eq:wass_equiv0}, where the \emph{pseudo-inverse} function for $\rho^{\mi,\mN}$ is given by
	\begin{equation*}
		X_{\rho^{\mi,\mN}}(m,t) = \sum_{k\in\mNN}\left( x^\mi_k(t) + \frac{m - k \sigma^\mi_N }{\rho^\mi_k(t)}  \right) \chi_{\big[k \sigma^\mi_N , (k+1)\sigma^\mi_N\big]}(m)\,.
	\end{equation*}
	Then, for any $t>s$, we have
	\begin{align*}
		d_{W^1}\big(\rho^{\mi,\mN}(t,\cdot),\rho^{\mi,\mN}(s,\cdot)\big) %&= \| \Gamma_{\rho^{\mi,N}}(m,t) - \Gamma_{\rho^{\mi,N}}(m,s)  \|_{L^1(0,1)}\\
		%&\le \sum_{k\in\mNN}\int \Bigg| x^\mi_k(t) - x^\mi_k(s) + (m - k \sigma^\mi_N) \bigg( \frac{1}{\rho^\mi_k(t)} - \frac{1}{\rho^\mi_k(s)} \bigg) \Bigg| \, \chi_{\big[k \sigma^\mi_N , (k+1)\sigma^\mi_N\big]}(m) \, dm \\
		&\le \sum_{k\in\mNN}\int_{k\sigma^\mi_N}^{(k+1)\sigma^\mi_N} \Big| x^\mi_k(t) - x^\mi_k(s) \Big|\, dm\\
		&\quad+ \sum_{k\in\mNN}\int_{k\sigma^\mi_N}^{(k+1)\sigma^\mi_N} \Bigg| (m - k \sigma^\mi_N) \bigg( \frac{1}{\rho^\mi_k(t)} - \frac{1}{\rho^\mi_k(s)} \bigg) \Bigg| \, dm \\
		%&\le \sum_{k\in\mNN}\sigma^\mi_N \Big| x^\mi_k(t) - x^\mi_k(s) \Big|+ \frac{(\sigma^\mi_N)^2}{2}  \Bigg| \bigg( \frac{1}{\rho^\mi_k(t)} - \frac{1}{\rho^\mi_k(s)} \bigg) \Bigg| \\
		&\le \sum_{k\in\mNN}\sigma^\mi_N \int_s^t \Big| \dot{x}^\mi_k(\tau) \Big|\, d\tau +\frac{(\sigma^\mi_N)^2}{2}  \int_s^t \Bigg|  \frac{d}{d\tau} \frac{1}{\rho^\mi_k(\tau)} \Bigg| \,d\tau\\
		%&\le \sum_{k\in\mNN}\sigma^\mi_N \int_s^t \Big| \frac{d}{d\tau}x^\mi_k(\tau) \Big|\, d\tau + \sum_{k\in\mNN}(\sigma^\mi_N)^2 \int_s^t \Bigg|  \frac{d}{d\tau} \frac{x^\mi_{k+1} - x^\mi_k}{\sigma^\mi_N}(\tau) \Bigg| \,d\tau\\
		&\le 3 \sigma^\mi_N \sum_{k\in\mNN} \int_s^t \Big| \dot{x}^\mi_k(\tau) \Big|\, d\tau\\
		&\leq C|t-s|,
	\end{align*}
	in view of Lemma \ref{lem:bound_velocity}. 
	
	Once proved the bounds described by equations \eqref{BVest} and \eqref{w_1bound}, we can apply Theorem \ref{Aubin Lions}, which concludes the proof.
\end{proof}

\begin{lemma}[Convergence of first momentum]\label{lem:mom}
	Given $\mu_{\rho^\mi}(t)$ and $\mu_{x^\mi}^{\mN}(t)$, respectively 
	\begin{equation*}
		\mu_{\rho^\mi}(t) = \frac{1}{\sigma^\mi|\Omega|}\int_\Omega y\,\rho^\mi(y,t)\,dy, \,\qquad \mu_{x^\mi}^{\mN}(t) = \frac{1}{N+1} \sum_{k=0}^{N} x^\mi_k(t),
	\end{equation*}
	we have that
	\begin{equation*}
		\lim_{N\to \infty} \big(\mu_{\rho^\mi}(t) - \mu_{x^\mi}^{\mN}(t)\big) = 0 
	\end{equation*}
	for all $t>0$.
\end{lemma}

\begin{proof}
	We recall that $x^i_0 + x^i_N = 0$, and $\Omega = [-1,1]$, from the definitions we have
	\begin{align*}
		\mu_{x^\mi}^{\mN}(t) & = \frac{1}{2(N+1)} \sum_{k\in\mNN}( x^\mi_{k+1} + x^\mi_k)= \frac{1}{2(N+1)} \sum_{k\in\mNN}\frac{(x^\mi_{k+1} + x^\mi_k) ( x^\mi_{k+1} - x^\mi_k)}{ x^\mi_{k+1} - x^\mi_k}\\
		%&= \frac{1}{2(N+1)} \sum_{k\in\mNN}\frac{(x^\mi_{k+1})^2 - (x^\mi_k)^2 }{ x^\mi_{k+1} - x^\mi_k}\\
		&= \frac{N}{(N+1)} \sum_{k\in\mNN}\frac{\rho^\mi_k }{\sigma^\mi} \frac{(x^\mi_{k+1})^2 - (x^\mi_k)^2}{2}= \frac{N}{(N+1)} \sum_{k\in\mNN} \frac{\rho^\mi_k(t)}{\sigma^\mi}\, \int_{-\infty}^{\infty} x\, \chi_{I^\mi_k} \,dx \\
		%&= \frac{N}{(N+1)} \int_{\Omega} \frac{x\,\rho^{\mi,\mN}(t,x)}{\sigma^\mi} \,dx \\
		&= \frac{N}{(N+1)}\, \mu_{\rho^{\mi,N}}(t) \,,
	\end{align*}
	where we used equation \eqref{eq:piec_discden}. We conclude that there exist a constant depending only on the domain $\Omega$ such that 
	\begin{align*}
		|\mu_{\rho^\mi}(t) - \mu_{x^\mi}^{\mN}(t) |  &\le \int_\Omega |\rho^\mi - \frac{N}{N+1}\rho^{\mi,N}|\, |x|\, dx\\
		&\le C(|\Omega|) \| \rho^\mi - \rho^{\mi,N} \|_{L^1} + o\left(\frac{C}{N}\right)\,,
	\end{align*}
	which concludes the proof.
\end{proof}  

We turn on the convergence of the approximated nodes.
\begin{prop}\label{prop:an}
	For any $T>0$ and $\mathbf{V}$ under assumption ($\mathbf{V}$). Then for any $\mi\in\mM$, there exists $a^\mi\in C([0,T]) $ such that $a_\mN^\mi\to a^\mi$ as $N\to\infty$ uniformly in $[0,T]$. Moreover, the limits $a^\mi$ satisfy \eqref{mild_a} for all $t\in[0,T]$.
\end{prop}
\begin{proof}
	We first notice that from \eqref{eq: evol agent} and the boundedness of $\mathbf{V}$ we have the uniform bound
	\begin{equation*}
		|a_\mN^\mi(t)|\leq |\bar{a}^\mi|+T\|\mathbf{V}\|_{\infty},\quad \mbox{ with }t\in[0,T],
	\end{equation*}
	for all $\mi\in\mM$ uniformly in $N$. Thus, there exist $a^\mi$ such that $a_\mN^\mi(t)$ punctually converges (up to subsequeces) to $a^\mi(t)$ as $N\to\infty$. Consider now $N_1,N_2\in\N$, then
	\begin{align*}
		\sum_{\mi\in\mM}|a_{\mN_1}^\mi(t)-a_{\mN_2}^\mi(t)|\leq & \sum_{\mi,\mj\in\mM} \int_0^t\left| \mathbf{V}(\mu_{x^\mi}^{\mN_1}, \mu_{x^\mj}^{\mN_1}, a_{\mN_1}^{\mij})-\mathbf{V}(\mu_{x^\mi}^{\mN_2}, \mu_{x^\mj}^{\mN_2}, a_{\mN_2}^{\mij})\right|d\tau\\
		\leq & C\sum_{\mi,\mj\in\mM}\int_0^t\left| \mu_{x^\mi}^{\mN_1}(\tau)-\mu_{x^\mi}^{\mN_2}(\tau)\right|+\left| \mu_{x^\mj}^{\mN_1}(\tau)-\mu_{x^\mj}^{\mN_2}(\tau)\right|\,d\tau\\
		& + C\sum_{\mj\in\mM}\int_0^t\left| a_{\mN_1}^{\mij}(\tau)-a_{\mN_2}^{\mij}(\tau)\right|\,d\tau\\
		% \leq & 2C\sum_{\mi\in\mM}\sum_{l=1,2}\| \rho^\mi-\rho^{\mi,\mN_l}\|_{L^1(\Omega_T)} + C\sum_{\mi,\mj\in\mM}\int_0^t\left| a_{\mN_1}^{\mij}(\tau)-a_{\mN_2}^{\mij}(\tau)\right|\,d\tau\\
		\leq & 2C\sum_{\mi\in\mM}\sum_{l=1,2}\| \rho^\mi-\rho^{\mi,\mN_l}\|_{L^1(\Omega_T)} + 2C\int_0^t\sum_{\mi\in\mM}|a_{\mN_1}^\mi(t)-a_{\mN_2}^\mi(t)|\,d\tau,
	\end{align*}
	where we used Lemma \ref{lem:mom} and the fact that by straightforward manipulations we have 
	\begin{align*}
		\sum_{\mi,\mj\in\mM} \left| a_{\mN_1}^{\mij}(\tau)-a_{\mN_2}^{\mij}(\tau)\right|=&% \sum_{\mi,\mj\in\mM}\left| \|a_{\mN_1}^{\mi}(\tau)-a_{\mN_1}^{\mj}(\tau)\|-\|a_{\mN_2}^{\mi}(\tau)-a_{\mN_2}^{\mj}(\tau)\|\right|\\
		%& =\sum_{\mi,\mj\in\mM}\left| \|a_{\mN_1}^{\mi}(\tau)\pm a_{\mN_2}^{\mi}(\tau)\pm a_{\mN_2}^{\mj}(\tau)-a_{\mN_1}^{\mj}(\tau)\|-\|a_{\mN_2}^{\mi}(\tau)-a_{\mN_2}^{\mj}(\tau)\|\right|\\
		% & \leq \sum_{\mi,\mj\in\mM}\left| \|a_{\mN_1}^{\mi}(\tau)-a_{\mN_2}^{\mi}(\tau)\|+\|a_{\mN_2}^{\mi}-a_{\mN_2}^{\mj}(\tau)\|\right.\\
		% &-\left.\|a_{\mN_2}^{\mi}(\tau)-a_{\mN_2}^{\mj}(\tau)\|+\|a_{\mN_2}^{\mj}(\tau)-a_{\mN_1}^{\mj}(\tau)\|\right|\\
		%& \leq \sum_{\mi,\mj\in\mM} \left|a_{\mN_1}^{\mi}(\tau)-a_{\mN_2}^{\mi}(\tau)\right|+\left|a_{\mN_2}^{\mj}(\tau)-a_{\mN_1}^{\mj}(\tau)\right|\\
		\leq 2\sum_{\mi\in\mM}|a_{\mN_1}^\mi(t)-a_{\mN_2}^\mi(t)|
	\end{align*}
	Thus, by Gronwall's type inequality we have
	%\begin{align*}
	%      u(t)\leq & \alpha + 2C\int_0^t u(\tau)\,d\tau\\
	%   \leq & \alpha e^{2Ct}
	%\end{align*}
	%we can deduce that
	\begin{align*}\sup_{t\in[0,T]}\sum_{\mi\in\mM}|a_{\mN_1}^\mi(t)-a_{\mN_2}^\mi(t)|\leq & 2C\sum_{\mi\in\mM}\sum_{l=1,2}\| \rho^\mi-\rho^{\mi,\mN_l}\|_{L^1(\Omega_T)} e^{2CT}.
	\end{align*}
	and then
	\begin{align*}
		\sup_{t\in[0,T]}|a_{\mN_1}^\mi(t)-a_{\mN_2}^\mi(t)|\leq & 2C\sum_{\mi\in\mM}\sum_{l=1,2}\| \rho^\mi-\rho^{\mi,\mN_l}\|_{L^1(\Omega_T)} e^{2CT},
	\end{align*}
	that ensure the uniform convergences of $a_{\mN_1}^\mi$ to $a^\mi,$ for all $\mi \in\mM$.
	
	In order to show that $a^\mi$ satisfies \eqref{mild_a} it is enough to observe that we can invoke the dominated convergence theorem since, by the continuity of $\mathbf{V}$ and the uniform converges proved we have that 
	\begin{equation*}
		\mathbf{V}(\mu_{x^\mi}^{\mN}, \mu_{x^\mj}^{\mN}, a_\mN^{\mij})\to \mathbf{V}(\mu_{\rho^\mi}, \mu_{\rho^\mj}, a^{\mij})\quad\mbox{a.e. in}  \quad t\in[0,T],
	\end{equation*} 
	and $\mathbf{V}(\mu_{x^\mi}^{\mN}, \mu_{x^\mj}^{\mN}, a_\mN^{\mij})$ is uniformly bounded w.r.t. $N$. Thus,
	\begin{equation*}
		\int_0^t\mathbf{V}(\mu_{x^\mi}^{\mN}, \mu_{x^\mj}^{\mN}, a_\mN^{\mij})\,d\tau \to \int_0^t \mathbf{V}(\mu_{\rho^\mi}, \mu_{\rho^\mj}, a^{\mij})\,d\tau,
	\end{equation*}
	for all $t\in[0,T]$.
\end{proof}
We now prove that the empirical measures associated to the solution of equation \eqref{eq: evol x} and the piece wise constant densities in equation \eqref{eq:piec_discden} share the same limit with respect to a suitable topology.
\begin{lemma}\label{lem:emph}
	For any $T>0$, the empirical measures associated to the solution of equation \eqref{eq: evol x}, defined by
	\begin{equation}\label{eq:empirical}
		\tilde{\rho}^{\mi,\mN}(t,x)=\sigma_N^\mi\sum_{k\in\mN}\delta_{x_k^\mi(t)}(x), \quad \mi\in\mM,
	\end{equation}
	satisfy
	\begin{equation*}
		d_{W^1}\left(\tilde{\rho}^{\mi,\mN}(t,\cdot),\rho^\mi(t,\cdot)\right)\to 0,\,\mbox{ as }\, N\to \infty,
	\end{equation*}
	for all $t\in[0,T]$, where $\rho^\mi$ is the limit obtained in Proposition \ref{comp}.
\end{lemma}
\begin{proof}
	We take again advantage from the isometry between the $1-$Wasserstain space for probability measures and the $L^1$ space in the space pseudo-inverse functions by noticing that the pseudo-inverse of an empirical measure is piece wise constant and then 
	\begin{align*}
		d_{W^1}\big(\tilde{\rho}^{\mi,\mN}(t,\cdot),\rho^{\mi,\mN}(t,\cdot)\big) = & \|X_{\tilde{\rho}^{\mi,\mN}}(t,\cdot)-X_{\rho^{\mi,\mN}}(t,\cdot)\|_{L^1([0,\sigma^\mi])}\\
		\leq & \sum_{k\in\mNN}\int_{k\sigma^\mi_N}^{(k+1)\sigma^\mi_N} \Bigg| (m - k \sigma^\mi_N)\frac{1}{\rho^\mi_k(t)} \Bigg| \, dm\\
		= & \frac{\sigma_N^\mi}{2}|\Omega|.
	\end{align*}
	The statement then follows from a triangulation argument.
\end{proof}

\subsection{Convergence to weak solutions}
With $\textbf{a}^\mN$ we denoted the vector $(a^\mi_\mN)_{\mi\in\mM}$, of agents with piece wise constant opinion distribution $\rho^{\mi,\mN}$, while $\textbf{a}$ is related to the vector of continuous distributions $\bm{\rho}$. This distinction is not stressed in the rest of the paper where the context does not allow misunderstanding. 
\begin{remark}
	We notice that evaluating the operator $\beta^\mi$ in \eqref{eq:empirical} we have
	\begin{equation}\label{eq:betaemph}
		\begin{aligned}
			\beta^\mi(\bm{\tilde{\rho}}^N,\textbf{a}^N;x_{k}^\mi) =& \sum_{\mj \in\mM} \int_\Omega \mathbf{A}^\mij(x_k^\mi,y,a_N^\mij)\tilde{\rho}^{\mj,N}(y,t)\,dy\\
			=& \sum_{\mj \in\mM} \sum_{l\in\mN}\sigma_N^\mj\mathbf{A}^\mij(x_k^\mi,x_l^\mj,a_N^\mij)\\
			=&\beta_k^\mi.
		\end{aligned}
	\end{equation}
	Moreover
	\begin{equation}\label{eq:betadev}
		\partial_x\beta^\mi\left(\tilde{\bm{\rho}}^{N},\textbf{a}^N;x\right)=\sum_{\mj\in\mM}\int_{\Omega}\partial_1 \bA^\mij(x,y;a_N^\mij)\tilde{\rho}^{\mj,N}(t,y)\,dy,
	\end{equation}
\end{remark}

\begin{lemma}\label{lem:conv}
	Let $T>0$, and consider the kernels $\bA^\mij$, $\bK^\mij$ under assumptions $(\bA)$ and $(\bK)$ respectively. Let $\rho^{\mi,\mN}$ and $\tilde{\rho}^{\mi,\mN}$ the sequences defined in equations \eqref{eq:piec_discden} and \eqref{eq:empirical} respectively, and their limits $\rho^\mi$ given by Proposition \ref{comp} and Lemma \ref{lem:emph}, for all $\mi \in \mM$. Then, for every $\zeta\in C_{0}^\infty\left(\Omega_T\right)$ we have
	\begin{equation}\label{eq:convtemp}
		\int_{\Omega_T}\rho^{\mi,\mN}\partial_t\zeta \,dx\,dt\to  \int_{\Omega_T}\rho^{\mi}\partial_t\zeta \,dx\,dt 
	\end{equation}
	\begin{equation}\label{eq:convdiff}
		\int_{\Omega_T}\Phi^\mi\left(\rho^{\mi,\mN}\right)\partial_x\beta^\mi\left(\tilde{\bm{\rho}}^{\mN},\textbf{a}^\mN;x\right)\partial_x\zeta \,dx\,dt\to \int_{\Omega_T}\Phi^\mi\left(\rho^{\mi}\right)\partial_x\beta^\mi\left(\bm{\rho},\textbf{a};x\right)\partial_x\zeta \,dx\,dt
	\end{equation}
	\begin{equation}\label{eq:convdiff2}
		\int_{\Omega_T}\Phi^\mi\left(\rho^{\mi,\mN}\right)\beta^\mi\left(\tilde{\bm{\rho}}^{\mN},\textbf{a}^\mN;x\right)\partial_{xx}\zeta \,dx\,dt\to \int_{\Omega_T}\Phi^\mi\left(\rho^{\mi}\right)\beta^\mi\left(\bm{\rho},\textbf{a};x\right)\partial_{xx}\zeta \,dx\,dt
	\end{equation}
	\begin{equation}\label{eq:convnon}
		\int_{\Omega_T}\rho^{\mi,\mN}\theta^\mi\left(\tilde{\bm{\rho}}^{\mN},\textbf{a}^\mN;x\right)\partial_x\zeta \,dx\,dt\to \int_{\Omega_T}\rho^{\mi}\theta^\mi\left(\bm{\rho},\textbf{a};x\right)\partial_x\zeta \,dx\,dt
	\end{equation}
	as $N\to\infty$.
\end{lemma}

\begin{proof}
	We only prove equation \eqref{eq:convdiff}, since equations \eqref{eq:convdiff2} and \eqref{eq:convnon} follow from similar argument, and equation \eqref{eq:convtemp} is a direct consequence of the $L^1$ strong compactness proved in Proposition \ref{comp}. We first split the terms as following
	\begin{align*}        &\left|\int_{\Omega_T}\left(\Phi^\mi\left(\rho^{\mi,\mN}\right)\partial_x\beta^\mi\left(\tilde{\bm{\rho}}^{\mN},\textbf{a}^\mN;x\right)-\Phi^\mi\left(\rho^{\mi}\right)\partial_x\beta^\mi\left(\bm{\rho},\textbf{a};x\right)\right)\partial_x\zeta(t,x)dx\,dt\right|\\
		&\leq \left|\int_{\Omega_T}\left(\Phi^\mi\left(\rho^{\mi,\mN}\right)-\Phi^\mi\left(\rho^{\mi}\right)\right)\partial_x\beta^\mi\left(\tilde{\bm{\rho}}^{\mN},\textbf{a}^\mN;x\right)\partial_x\zeta(t,x)\,dx\,dt \right|\\&\quad+\left|\int_{\Omega_T}\Phi^\mi\left(\rho^{\mi}\right)\left(\partial_x\beta^\mi\left(\tilde{\bm{\rho}}^{\mN},\textbf{a}^\mN;x\right)-\partial_x\beta^\mi\left(\bm{\rho},\textbf{a};x\right)\right)\partial_x\zeta(t,x)dx\,dt\right|\\
		& =|I|+|II|.
	\end{align*}
	We now treat the two terms separately. 
	Assumption $(\bA)$ and equation \eqref{eq:betadev} ensure the following bound
	\begin{align*}
		|I|\leq & \sum_{\mj\in\mM} \int_{\Omega_T}\left|\Phi^\mi\left(\rho^{\mi,\mN}\right)-\Phi^\mi\left(\rho^{\mi}\right)\right|\int_{\Omega}\left|\partial_1\bA^\mij(x,y;a^\mij)\tilde{\rho}^{\mj,\mN}(t,y)\,dy\right|\left|\partial_x\zeta(t,x)\right|\,dx\,dt\\
		\leq & C\|\rho^{\mi,\mN}-\rho^{\mi}\|_{L^1(\Omega_T)}.
	\end{align*}
	where $C$ is a constant depending on $\left\|\partial_1\bA^\mij\right\|_{\infty}$, $\left\|\partial_x\zeta\right\|_{\infty}$, $Lip\left(\Phi^\mi\right)$, and $\sigma^\mM$. In order to bound the second integral let us introduce $\Pi^{\mi,\mN}$ an optimal transport plan between $\tilde{\rho}^{\mi,\mN}$ and $\rho^\mi$. Then we have
	\begin{align*}
		|II|\leq & \sum_{\mj\in\mM} \int_{\Omega_T}\Phi^\mi\left(\rho^{\mi}\right)\int_{\Omega^2}\left|\partial_1\bA^\mij(x,y;a_\mN^\mij)-\partial_1\bA^\mij(x,z;a^\mij)\right|\,d\Pi ^{\mi,\mN}(y,z)|\partial_x\zeta|\,dx\,dt\\
		\leq & C\sum_{\mj\in\mM} \int_{\Omega_T}\int_{\Omega\times\Omega}\left|y-z\right|\,d\Pi ^{\mi,\mN}(y,z)\,dx\,dt\\ \leq  & C M|\Omega|\int_0^T d_{W^1}(\tilde{\rho}^{\mi,\mN}(t,\cdot),\rho^\mi(t,\cdot))\,dt,
	\end{align*}
	where in this case $C$ is a constant depending on $\|\rho^\mi\|_\infty$, $\left\|\partial_x\zeta\right\|_{\infty}$, and the constant $c_{1,A}$ from Assumption $(\bA)$. The convergences in Proposition \ref{comp} and Lemma \ref{lem:emph} ensure that equation \eqref{eq:convdiff} holds. 
\end{proof}

We are now in the position of proving that the limit densities and nodes satisfy the weak formulation in the sense of Definition \ref{def:weak}. More precisely we are going to show that for $N\to+\infty$ we have
\begin{equation}\label{weaklimit}
	\begin{aligned}
		\int_{\Omega_T}&\rho^{\mi,\mN}\partial_t\zeta+\rho^{\mi,\mN}\theta^\mi(\bm{\tilde{\rho}}^\mN,\textbf{a}^\mN;x) \partial_x \zeta\, dx\,dt\\
		&+\int_{\Omega_T}\Phi^\mi\left(\rho^{\mi,\mN}\right) \left( \partial_x\beta^\mi(\bm{\tilde{\rho}}^\mN,\textbf{a}^\mN;x)   \partial_x \zeta+\beta^\mi(\bm{\tilde{\rho}}^\mN,\textbf{a}^\mN;x)   \partial_{xx} \zeta\right)\, dx\,dt\to 0        
	\end{aligned}
\end{equation}
that combined with the convergences in Lemma \ref{lem:conv} gives the assertion. We state the following 
\begin{prop}
	Given $M\in\N$, and $T>0$ fixed, for all $\mi,\mj\in \mM$ consider $\bA^{\mij}$, $\bK^{\mij}$, $\mathbf{V}$, and $\Phi^\mi$ under assumptions $(\bA)$, $(\bK)$, $(\textbf{V})$ and $(Dif)$ respectively. Let $\bar{\rho}^\mi:\Omega \to \R$ under assumptions $(In1)$ and $(In2)$ for all $\mi\in\mM$. \\ Then, for all $\mi \in \mM$ the densities $\rho^{\mi,\mN}$ introduced in equation \eqref{eq:piec_discden} and $\hat{\rho}^{\mi,\mN}$ introduced in equation \eqref{eq:empirical} satisfy equation \eqref{weaklimit} as $N\to\infty$, for all $\zeta\in C_0^{\infty}(\Omega_T)$. 
\end{prop}
\begin{proof}
	We start considering the term involving the time derivative. By definition of $\rho^{\mi,\mN}$ in equation \eqref{eq:piec_discden}, a discrete integration by parts and Fundamental Theorem of Calculus  give
	\begin{align*}  \int_{\Omega_T}\rho^{\mi,\mN}\partial_t\zeta=&\sum_{k\in\mNN}\int_0^T\rho_k^\mi(t)\int_{I_k^\mi(t)}\partial_t\zeta(t,x)\,dx\,dt\\
		=&\sum_{k\in\mNN}\int_0^T\rho_k^\mi(t)\dot{x}_{k+1}^\mi(t)\left(\int_{I_k^\mi(t)}\zeta(t,x)\,dx-\zeta(t,x_{k+1}^\mi(t))\right)\,dt\\
		&-\sum_{k\in\mNN}\int_0^T\rho_k^\mi(t)\dot{x}_{k}^\mi(t)\left( \int_{I_k^\mi(t)}\zeta(t,x)\,dx-\zeta(t,x_{k}^\mi(t))\right)\,dt.
	\end{align*}
	A second order expansion of $\zeta$ around $x_{k+1}^\mi(t)$ in the first average integral and around $x_{k}^\mi(t)$ in the second average integral produces
	\begin{align*}
		\int_{\Omega_T}\rho^{\mi,\mN}\partial_t\zeta=&-\frac{\sigma_N^\mi}{2}\sum_{k\in\mNN}\int_0^T\dot{x}_{k+1}^\mi(t)\partial_x\zeta(t,x_{k+1}^\mi)\,dt\\
		&+\frac{1}{2}\sum_{k\in\mNN}\int_0^T\rho_k^\mi(t)\dot{x}_{k+1}^\mi(t)\int_{I_k^\mi(t)}\partial_{xx}\zeta(t,\hat{x}_{k+1}^\mi)\left(x-x_{k+1}^\mi(t))\right)^2\,dx\,dt\\
		&-\frac{\sigma_N^\mi}{2}\sum_{k\in\mNN}\int_0^T\dot{x}_{k}^\mi(t)\partial_x\zeta(t,x_{k}^\mi)\,dt\\
		&-\frac{1}{2}\sum_{k\in\mNN}\int_0^T\rho_k^\mi(t)\dot{x}_{k}^\mi(t)\int_{I_k^\mi(t)}\partial_{xx}\zeta(t,\hat{x}_{k}^\mi)\left(x-x_{k}^\mi(t))\right)^2\,dx\,dt,
	\end{align*}
	where $\hat{x}_{k}^\mi$ and $\hat{x}_{k+1}^\mi$ are points in $\left[x_{k}^\mi,x\right]$ and equation $\left[x,x_{k+1}^\mi\right]$ respectively. We now combine the first and third term on the r.h.s.~above and use \eqref{eq:DPA} in order to obtain
	\begin{equation*}
		-\frac{\sigma_N^\mi}{2}\sum_{k\in\mNN}\int_0^T\left(\dot{x}_{k}^\mi(t)\partial_x\zeta(t,x_{k}^\mi)+\dot{x}_{k+1}^\mi(t)\partial_x\zeta(t,x_{k+1}^\mi)\right)\,dt = A_1+A_2,
	\end{equation*}
	where
	\begin{align*}
		A_1 = - \frac{1}{2}\sum_{k\in\mNN}\int_0^T&\left(\beta^\mi_k\left(\Phi^\mi(\rho^\mi_{k-1}) - \Phi^\mi(\rho^\mi_{k}) \right) \partial_x\zeta(t,x_{k}^\mi)\right.\\
		&+\left.\beta^\mi_{k+1} \left(\Phi^\mi(\rho^\mi_{k}) - \Phi^\mi(\rho^\mi_{k+1}) \right) \partial_x\zeta(t,x_{k+1}^\mi)\right)\,dt,
	\end{align*}
	and
	\begin{equation*}
		A_2=  -  \frac{\sigma_N^\mi}{2}\sum_{k\in\mNN}\int_0^T\left(\theta_{k}^\mi(t)\partial_x\zeta(t,x_{k}^\mi)+\theta_{k+1}^\mi(t)\partial_x\zeta(t,x_{k+1}^\mi)\right)\,dt.
	\end{equation*}
	We now combine the integral $A_1$ with the two terms involving the diffusion in equation \eqref{weaklimit} in order to show that
	\begin{equation*}
		A_1 + \int_{\Omega_T}\Phi^\mi\left(\rho^{\mi,\mN}\right)\partial_x \left( \beta^\mi(\bm{\tilde{\rho}}^\mN,\textbf{a}^\mN;x)   \partial_x \zeta (t,x)\right)\, dx\,dt=0.
	\end{equation*}
	Invoking equation \eqref{eq:betaemph} and the fact that $\partial_x\zeta(t,x_{0}^\mi)=\partial_x\zeta(t,x_{N}^\mi)=0$ we can compute
	\begin{align*}
		&\int_{\Omega_T}\Phi^\mi\left(\rho^{\mi,\mN}\right)\partial_x \left( \beta^\mi(\bm{\tilde{\rho}}^\mN,\textbf{a}^\mN;x)   \partial_x \zeta (t,x)\right)\, dx\,dt\\
		=&\sum_{k\in\mNN}\int_0^T\Phi^\mi\left(\rho_k^{\mi}\right) \int_{I_k^\mi(t)}\partial_x \left( \beta^\mi(\bm{\tilde{\rho}}^\mN,\textbf{a}^\mN;x)   \partial_x \zeta(t,x)\right)\, dx\,dt\\
		=&\sum_{k\in\mNN}\int_0^T\Phi^\mi\left(\rho_k^{\mi}\right) \left( \beta^\mi(\bm{\tilde{\rho}}^\mN,\textbf{a}^\mN;x_{k+1}^\mi)   \partial_x \zeta(t,x_{k+1}^\mi)-\beta^\mi(\bm{\tilde{\rho}}^\mN,\textbf{a}^\mN;x_{k}^\mi)   \partial_x \zeta(t,x_{k}^\mi)\right)\,dt\\
		=&\sum_{k\in\mNN}\int_0^T\Phi^\mi\left(\rho_k^{\mi}\right) \left( \beta_{k+1}^\mi\partial_x \zeta(t,x_{k+1}^\mi)-\beta_k^\mi \partial_x \zeta(t,x_{k}^\mi)\right)\,dt\\
		=&\sum_{k\in\mNZN}\int_0^T\beta^\mi_k \left(\Phi^\mi(\rho^\mi_{k-1}) - \Phi^\mi(\rho^\mi_{k}) \right) \partial_x\zeta(t,x_{k}^\mi)\,dt,
	\end{align*}
	which can be combined with $A_1$ by shifting the indexes and using again the fact that the test function vanishes at the boundary.

	Recalling the definition of $\theta_k^\mi$ in equation \eqref{eq: def theta}, rearranging the indexes in $A_2$, using the fact that $\partial_x\zeta$ vanishes at the boundary of $\Omega$ and summing with the second term in equation \eqref{weaklimit} we obtain 
	\begin{align*}
		&\sum_{k\in\mNZN}\sum_{\mi\in\mM}\sum_{l\in\mNN}\sigma_N^\mi\sigma_N^\mj\int_0^T\bK^{\mij}(x_k^\mi,x_l^\mj;a_N^\mij)\partial_x\zeta(t,x_{k}^\mi)\,dt\\&+ \sum_{k\in\mNN}\sum_{\mi\in\mM}\sum_{l\in\mNN}\sigma_N^\mj\int_0^T\rho_k^{\mi}(t) \int_{I_{k}^\mi} \bK^\mij(x,x_l^\mj;a_N^\mij) \partial_x \zeta(t,x)\, dx\,dt. 
	\end{align*}
	A first order expansion on $\partial_x\zeta$ around $x_k^\mi$ for $\hat{\hat{x}}_k^\mi\in \left[x_k^\mi,x\right]$, together with the definition of $\rho_k^\mi$ and assumption $(\bK)$ yield
	\begin{align*}
		&\left|\sum_{k\in\mNZN}\sum_{\mi\in\mM}\sum_{l\in\mNN}\sigma_N^\mj\int_0^T\rho_k^\mi\partial_x\zeta(t,x_{k}^\mi)\int_{I_{k}^\mi} \bK^{\mij}(x_k^\mi,x_l^\mj;a_N^\mij)-\bK^\mij(x,x_l^\mj;a_N^\mij)\,dx\,dt\right.\\&\left.+ \sum_{k\in\mNN}\sum_{\mi\in\mM}\sum_{l\in\mNN}\sigma_N^\mj\int_0^T\rho_k^{\mi}(t) \int_{I_{k}^\mi} \bK^\mij(x,x_l^\mj;a_N^\mij) (x-x_k^\mi)\partial_{xx} \zeta(t,\hat{\hat{x}}_k^\mi)\, dx\,dt\right|\\
		&\leq \sum_{k\in\mNZN}\sum_{\mi\in\mM}\sum_{l\in\mNN}\sigma_N^\mj \left(c_\bK\|\partial_x\zeta\|_\infty+\|\bK^\mij\|_\infty \|\partial_{xx}\zeta\|_\infty\right)\int_0^T\rho_k^\mi\int_{I_{k}^\mi}|x_k^\mi-x|\,dx\,dt\\
		&=\sum_{k\in\mNZN}\sum_{\mi\in\mM}\frac{\sigma^\mj}{2} \left(c_\bK\|\partial_x\zeta\|_\infty+\|\bK^\mij\|_\infty \|\partial_{xx}\zeta\|_\infty\right)\int_0^T\sigma_N^\mi|I_k^\mi|(t)\,dt\\
		&\leq \sigma_N^\mi\sum_{\mi\in\mM}\frac{\sigma^\mj}{2} \left(c_\bK\|\partial_x\zeta\|_\infty+\|\bK^\mij\|_\infty \|\partial_{xx}\zeta\|_\infty\right)T|\Omega|,
	\end{align*}
	that vanishes as $N\to\infty$ toghether with $\sigma_N^\mi$.
	We are now left in showing that the remainder, i.e.~
	\begin{align*}
		R_{k,k+1}^\mi=&\frac{1}{2}\sum_{k\in\mNN}\int_0^T\rho_k^\mi(t)\dot{x}_{k+1}^\mi(t)\int_{I_k^\mi(t)}\partial_{xx}\zeta(t,\hat{x}_{k+1}^\mi)\left(x-x_{k+1}^\mi(t))\right)^2\,dx\,dt\\
		&-\frac{1}{2}\sum_{k\in\mNN}\int_0^T\rho_k^\mi(t)\dot{x}_{k}^\mi(t)\int_{I_k^\mi(t)}\partial_{xx}\zeta(t,\hat{x}_{k}^\mi)\left(x-x_{k}^\mi(t))\right)^2\,dx\,dt,
	\end{align*}goes to zero. We first notice that for all $h,k\in\mNN$ we have
	\begin{equation*}
		\int_{I_k^\mi(t)}\partial_{xx}\zeta(t,\hat{x}_{h}^\mi)\left(x-x_{k}^\mi(t))\right)^2\,dx\leq \|\partial_{xx}\zeta\|_\infty \frac{|I_k^\mi|^3(t)}{3}\leq C
	\end{equation*}
	for some constant $C>0$, then
	\begin{align*}
		|R_{k,k+1}^\mi|\leq &\frac{\|\partial_{xx}\zeta\|_\infty}{2}\sum_{k\in\mNN}\int_0^T\rho_k^\mi(t)\left(|\dot{x}_{k+1}^\mi(t)|+|\dot{x}_{k}^\mi(t)|\right)\frac{|I_k^\mi|^3(t)}{3}\,dt\\
		\leq & C\|\partial_{xx}\zeta\|_\infty\sum_{k\in\mNN}\int_0^T\rho_k^\mi(t)|I_k^\mi|^3(t)\,dt\\
		=& \sigma_N^\mi C\|\partial_{xx}\zeta\|_\infty|\Omega|T,
	\end{align*}
	where the second inequality holds in view of Lemma \ref{lem:bound_velocity} and the BV bound on $\rho_k^\mi$.
\end{proof}

\section{Modelling and simulation}\label{sec: modelling}
As mentioned in the introduction, the aim of this model is to describe a few aspects typical of the interactions on social networks and social media. In particular, those aspects related to the evolution of the network and those concerning the homophilia and heterophobia. The goal is to observe if we get polarization of the opinions and fragmentation of the population only modelling the processes ruling the rewiring of the network, and the epistemic process. We do not model the polarization itself, we stick to the description of the opinion formation dynamics that drive the works on formation of echo chambers and epistemic bubbles. In our model we introduce two main concepts: the \textbf{attitude areas}, the \textbf{euclidean network}. 

\subsection{Attitude areas}
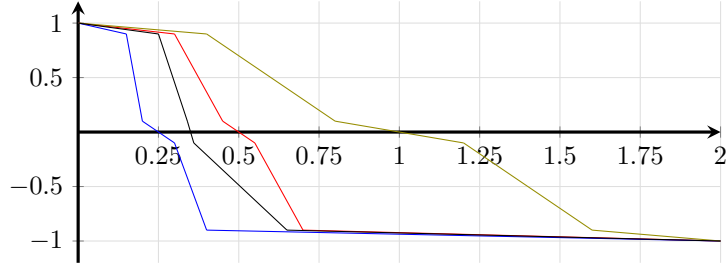
\begin{figure}
	\centering
	\begin{tikzpicture} 
		\begin{axis}[
			axis line style = very thick,
			axis lines=middle,
			xmin = 0, xmax = 2,
			ymin = -1.2, ymax = 1.2,
			xtick distance = 0.25,
			ytick distance = 0.5,
			grid = major,
			major grid style = {lightgray!50},
			width = 0.8\textwidth,
			height = 0.4\textwidth]

			\addplot[red,
			domain = 0:0.3,
			] {1 - 0.1 * x/(0.3) };
			\addplot[red,
			domain = 0.3:0.45,
			] {0.1 + 0.8 * (1 - (x-0.3)/(0.45-0.3)};
			\addplot[red,
			domain = 0.45:0.55,
			] {-0.1 + 0.2 * (1 - (x - 0.45)/(0.55-0.45))};
			\addplot[red,
			domain = 0.55:0.7,
			] {-0.9 + 0.8 * (1 - (x - 0.55)/(0.7-0.55))};
			\addplot[red,
			domain = 0.7:2,
			] {-0.9 - 0.1 * (x - 0.7)/(2-0.7)};

			\addplot[blue,
			domain = 0:0.15,
			] {1 - 0.1 * x/(0.15) };
			\addplot[blue,
			domain = 0.15:0.2,
			] {0.1 + 0.8 * (1 - (x-0.15)/(0.2-0.15)};
			\addplot[blue,
			domain = 0.2:0.3,
			] {-0.1 + 0.2 * (1 - (x - 0.2)/(0.3-0.2))};
			\addplot[blue,
			domain = 0.3:0.4,
			] {-0.9 + 0.8 * (1 - (x - 0.3)/(0.4-0.3))};
			\addplot[blue,
			domain = 0.4:2,
			] {-0.9 - 0.1 * (x - 0.4)/(2-0.4)};

			\addplot[olive,
			domain = 0:0.4,
			] {1 - 0.1 * x/(0.4) };
			\addplot[olive,
			domain = 0.4:0.8,
			] {0.1 + 0.8 * (1 - (x-0.4)/(0.8-0.4)};
			\addplot[olive,
			domain = 0.8:1.2,
			] {-0.1 + 0.2 * (1 - (x - 0.8)/(1.2-0.8))};
			\addplot[olive,
			domain = 1.2:1.6,
			] {-0.9 + 0.8 * (1 - (x - 1.2)/(1.6 - 1.2))};
			\addplot[olive,
			domain = 1.6:2,
			] {-0.9 - 0.1 * (x - 1.6)/(2-1.6)};

			\addplot[
			domain = 0:0.25,
			] {1 - 0.1 * x/(0.25) };
			\addplot[
			domain = 0.25:0.34,
			] {0.1 + 0.8 * (1 - (x-0.25)/(0.34-0.25)};
			\addplot[
			domain = 0.34:0.36,
			] {-0.1 + 0.2 * (1 - (x - 0.34)/(0.36-0.34))};
			\addplot[
			domain = 0.36:0.65,
			] {-0.9 + 0.8 * (1 - (x - 0.36)/(0.65-0.36))};
			\addplot[
			domain = 0.65:2,
			] {-0.9 - 0.1 * (x - 0.65)/(2-0.65)};
			
		\end{axis}
	\end{tikzpicture}
	\caption{Attraction/repulsion function} \label{fig: scala}
	\medskip
	\small
	The positive values of the function coincide with the attraction, on the other hand while the function has negative values it describes repulsion between the agents' opinion, which could bring to radicalization or polarization. Due to the choice of the domain we have that $s\in[0,2]$. The different colors coincide with the following definitions of the attitude intervals:\\
	\color{blue}Blue\color{black}: $r_\textit{f}=0.15$, $r_\textit{a}=0.20$, $r_\textit{r}=0.30$, $r_\textit{l}=0.40$.\\ Black: $r_\textit{f}=0.25$, $r_\textit{a}=0.34$, $r_\textit{r}=0.36$, $r_\textit{l}=0.65$.\\ \color{red}Red\color{black}: $r_\textit{f}=0.30$, $r_\textit{a}=0.45$, $r_\textit{r}=0.55$, $r_\textit{l}=0.70$.\\ \color{olive}Olive\color{black}: $r_\textit{f}=0.40$, $r_\textit{a}=0.80$, $r_\textit{r}=1.20$, $r_\textit{l}=1.60$.
\end{figure}

In this section we model and simulate a possible choice of opinion dynamic. In particular we focus our attention on a model that takes into account the 
heterophobia and homophilia dynamics. The model that we propose is based on  five attitude areas: \textit{attraction}/\textit{homophilia}, \textit{curiosity}, \textit{indifference}, \textit{mistrust}, \textit{repulsion}/\textit{heterophobia}. When two agents interact, their attitude depends on the distance between their respective opinions. We do not consider the attitude depending locally on the opinion itself, this could be a possible improvement that describes the low propensity of changing extreme beliefs. The interaction can be attractive, i.e.~the agents moves toward a position of consensus looking for a compromise, or can be repulsive, i.e.~each agent changes its own opinion moving farther from the one of the other agent.
We consider the DPA structure of equations (\ref{eq: evol x}, \ref{eq: evol agent}, \ref{eq: def beta}, \ref{eq: def theta}). The operator describing this phenomenon has the following structure
\begin{equation}\label{func: int kernel reaction 1}
	\mathbf{K}^\mij(w,v,a^\mij) = \omega(a^\mij) \, \zeta(\mu^\mi - \mu^\mj)\, (v - w)
\end{equation}
where the function $\omega$ depends on the network connections, and $\zeta$ is the \textit{attitude function} that measures the distance between the agents' mean opinion. Being $s$ the distance between the mean opinions, the five attitude areas coincide with the following intervals
$$  \begin{aligned}
	&s \in (0,r_{\textit{f}})  \quad &&\text{strong attraction, homophilia} &&(r_\textit{friends})\\
	&s \in (r_{\textit{f}},r_{\textit{a}}) \quad &&\text{curiosity} &&(r_\textit{attraction})\\    
	&s \in (r_{\textit{a}},r_{\textit{r}}) \quad &&\text{indifference} \\
	&s \in (r_{\textit{r}},r_{\textit{l}}) \quad &&\text{mistrust} &&(r_\textit{repulsion})\\
	&s \in (r_{\textit{l}},|\Omega|) \quad &&\text{repulsion, heterophobia}&&(r_\textit{limit})\,.
\end{aligned}
$$
The function $\zeta$ is given by
$$
\zeta(s) = 
\begin{cases}
	\begin{aligned}
		&1 - \frac{1}{10} \frac{|s|}{r_{\textit{f}}} \quad &if \quad & |s| < r_{\textit{f}}\\
		&0.1 + \frac{8}{10} \bigg[1 - \frac{|s| - r_{\textit{f}}}{r_{\textit{a}} - r_{\textit{f}}} \bigg] \quad &if \quad & r_{\textit{f}} \le |s| < r_{\textit{a}}\\    
		& -0.1 + \frac{2}{10} \bigg[1 - \frac{|s| - r_{\textit{a}}}{r_{\textit{r}} - r_{\textit{a}}} \bigg] \quad &if \quad & r_{\textit{a}} \le |s| < r_{\textit{r}}\\
		&-0.9 + \frac{8}{10} \bigg[1 - \frac{|s| - r_{\textit{r}}}{r_{\textit{l}} - r_{\textit{r}}} \bigg] \quad &if \quad & r_{\textit{r}} \le |s| < r_{\textit{l}}\\
		&-0.9 - \frac{1}{10} (|s| - r_{\textit{l}})  \quad &if \quad & r_{\textit{l}} \le |s| \,,
	\end{aligned}
\end{cases}
$$
depending on the values of the extreme of the intervals it looks like those in Figure \ref{fig: scala}.

\subsection{Diffusion}
Moreover, the opinion dynamic is not ruled only by the direct interaction with the connected agents in the network. We continuously get inputs from all the media, this phenomenon is described in this model by the following operator
\begin{equation}\label{func: int kernel reaction 2}
	\mathbf{A}^\mij(w,v,a^\mij) = |\mu^\mj - w|^2 \,.
\end{equation}
In this case the interaction is not filtered by the network connections, i.e.~by $\omega$, and it is not affected by the attitude. This operator plays the role of the diffusion mobility, the opinion tends to diffuse the more the inputs are far from it, and it is not affected by the diffusion when the inputs coincides with the opinion itself. The underlying idea is that we cannot process actively - i.e.~through $\zeta$ - all the information that we get. Those inputs that we cannot elaborate they influence our opinion distribution smoothing it depending on the distance between the input and our mean opinion.

\subsection{Evolution of the network}
\begin{figure}
	\centering
	\includegraphics[width=0.6\textwidth]{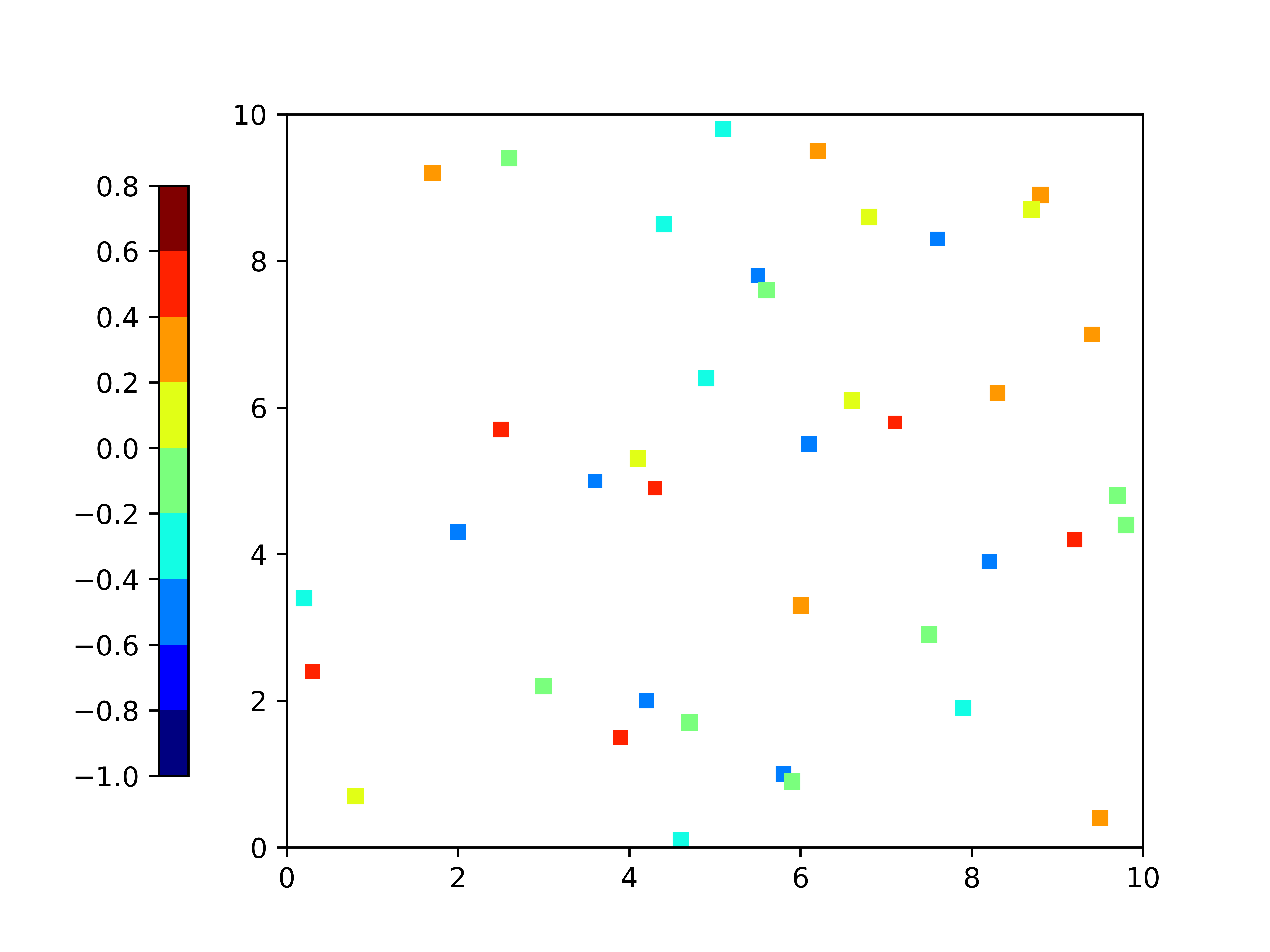}
	\caption{Initial network condition}\label{fig: ic network}
	\medskip
	\small
	In this figure the initial agents coordinates belong to $[0,10]^2$. The colors describe the mean opinion of each agent, which belongs to the interval $[-1,1]$. The number of agents is $N=40$. The agents coordinates are uniformly random distributed on each axes. The dimension of each square is proportional to the total mass $\sigma$ of each agent, in this case they all almost coincide.
\end{figure}

The interaction on the network is ruled by the distance between the agents in the network space, and by the distance between their respective mean opinions. The radius $r_{\textit{loc}}$ is the discriminant of the local interaction,
\begin{equation}
	\omega(a^\mij) = \begin{cases}
		1 \quad if \quad a^\mij \le \rho_{\textit{loc}}\\
		0 \quad if \quad a^\mij > \rho_{\textit{loc}}
	\end{cases} \,.
\end{equation}
The agents $a^\mi$ and $a^\mj$ interact if their distance on the network is smaller or equal to $r_{\textit{loc}}$.
The global interaction coincides with the radius $\rho_{\textit{loc}}=+\infty$.
The initial condition used for our simulation is the one given in Figure \ref{fig: ic network}.

\begin{figure}
	\centering
	\subfloat[$\rho_\textit{loc}=3.5$]{\includegraphics[width = 0.33\textwidth]{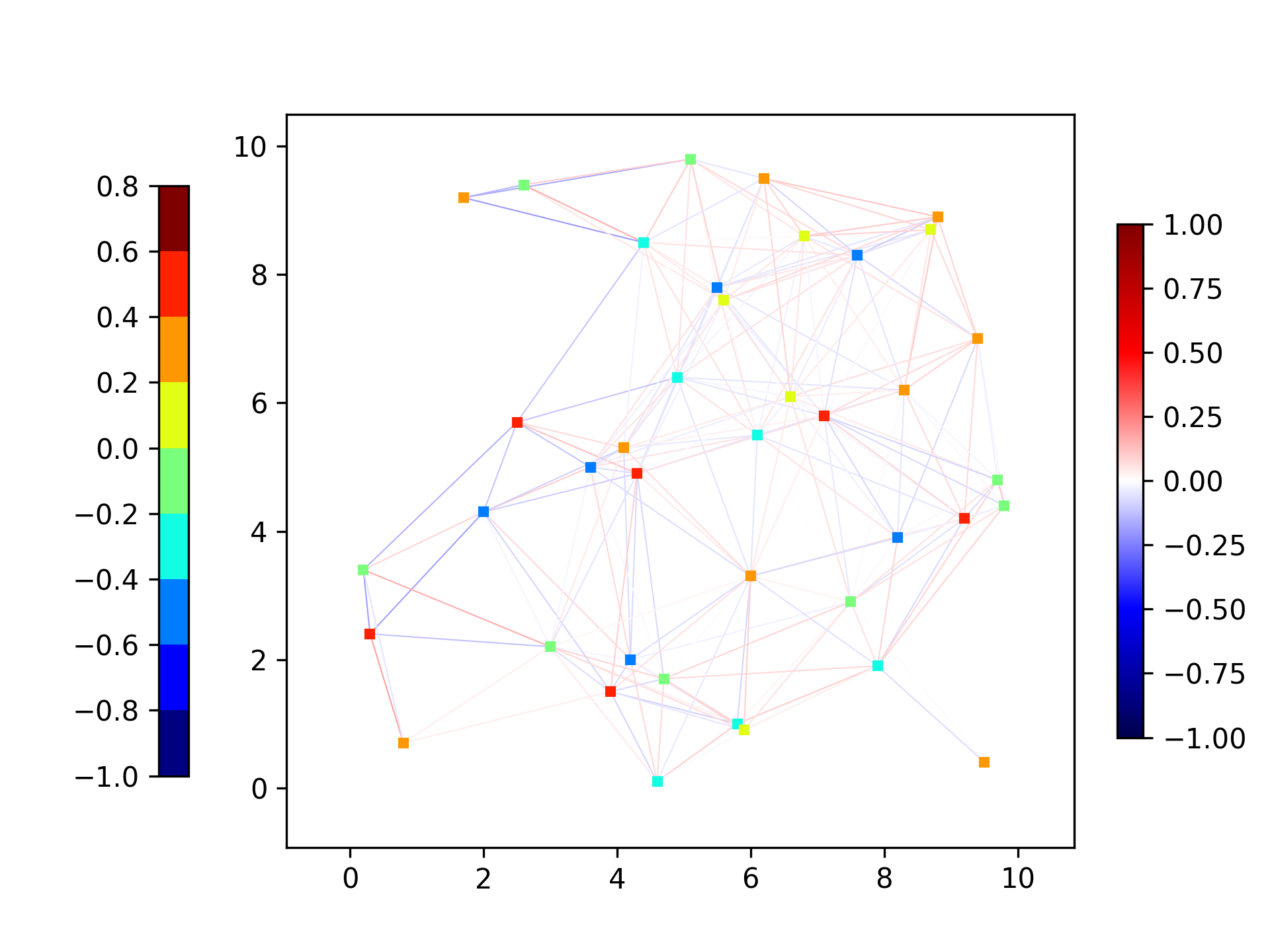}}
	\subfloat[$\rho_\textit{loc}=5$]{\includegraphics[width = 0.33\textwidth]{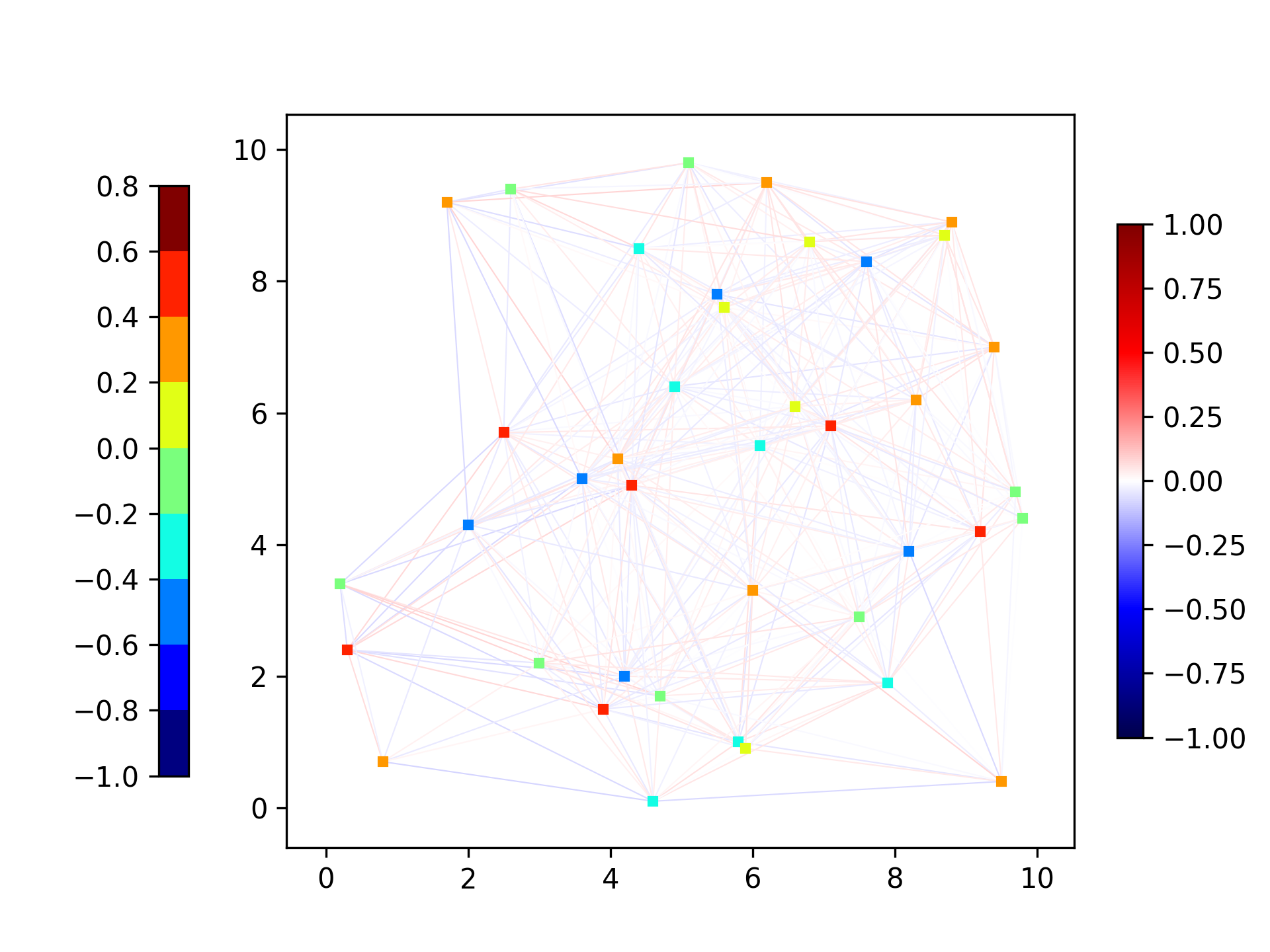}}
	\subfloat[$\rho_\textit{loc}=10$]{\includegraphics[width = 0.33\textwidth]{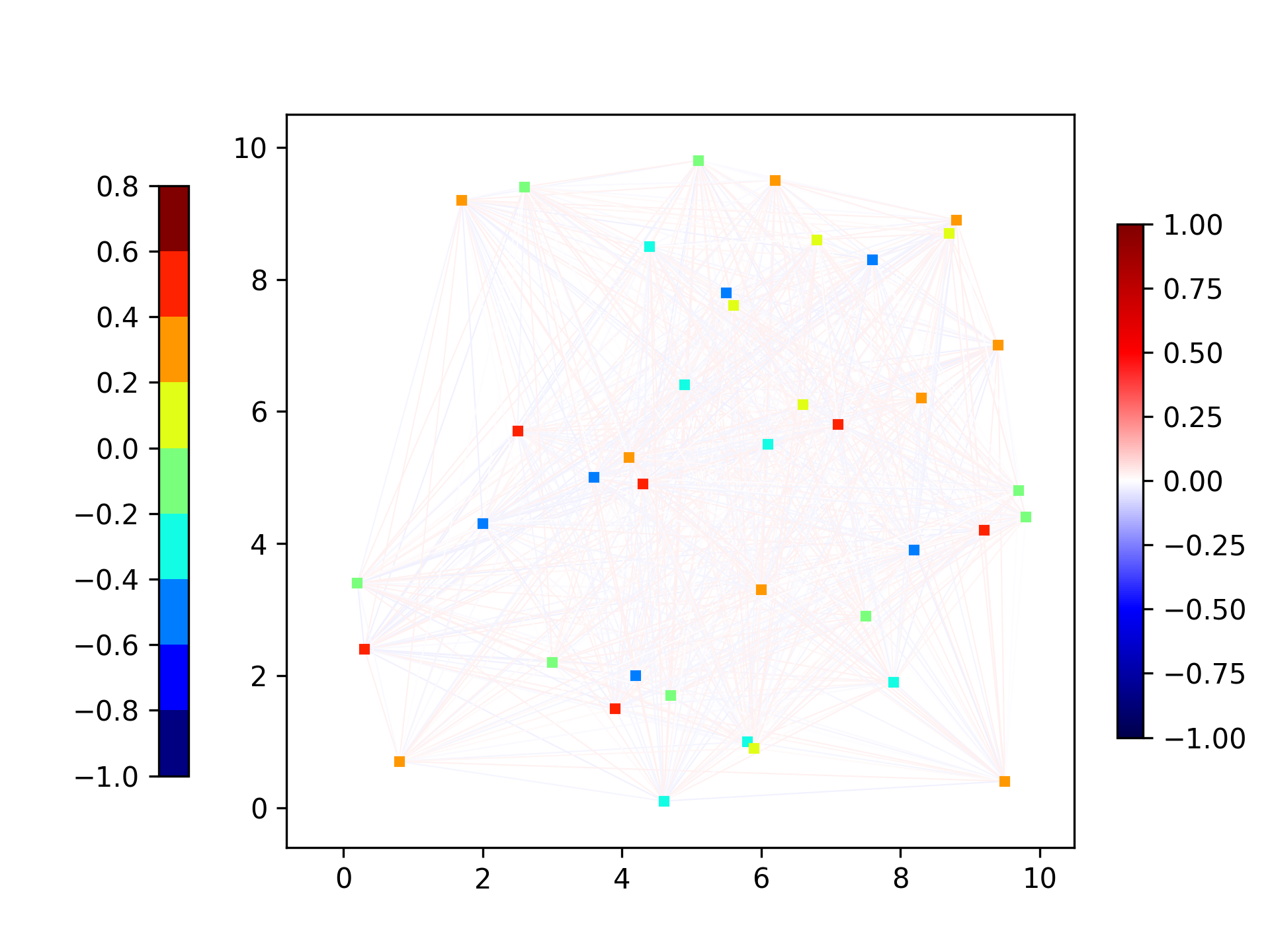}}
	\caption{Local initial network interaction}
	\label{fig: ic network connection}
	\medskip
	\small
	The agents interact if they are connected by a link. The magnitude and the sign of the connection ranges in $[-1,1]$ and are described by the legend on the right of the pictures. In this case the attitude areas are given by the parameters : $r_\textit{f}=0.25$, $r_\textit{a}=0.34$, $r_\textit{r}=0.36$, $r_\textit{l}=0.65$, i.e.~the black function in Figure \ref{fig: scala}.
\end{figure}

Considering 3 different radius of interaction the network connections change as in Figure \ref{fig: ic network connection}. 

The evolution of the agents is ruled by the following operator
\begin{equation}\label{func: int kernel network 1}
	\mathbf{V}(\mu_{\rho^\mi}, \mu_{\rho^\mj}, a^\mij) = \sum_{l=1}^{2} \sum_j \zeta|\mu_{\rho^\mi} - \mu_{\rho^\mj}| \, \omega(a^\mij) \, (a^\mi - a^\mj) \mathbf{e}_l\,,
\end{equation}
with $\mathbf{e}_l$ being the $j$-th versor.

\subsection{Radicalization, polarization, and fragmentation}
We call radicalization the tendency of the opinions to cluster on a value that does not coincide with the global consensus. While the polarization implies - on top of the radicalization - that the distributions move towards the extreme opinions, in our case $\pm 1$.\\
We consider the attitude areas given by the following intervals: $r_\textit{f}=0.25$, $r_\textit{a}=0.34$, $r_\textit{r}=0.36$, $r_\textit{l}=0.65$. This setting coincides with the $\zeta$ function described by the black line in Figure 
\ref{fig: scala}.

\begin{figure}
	\centering
	\includegraphics[width = 0.6\textwidth]{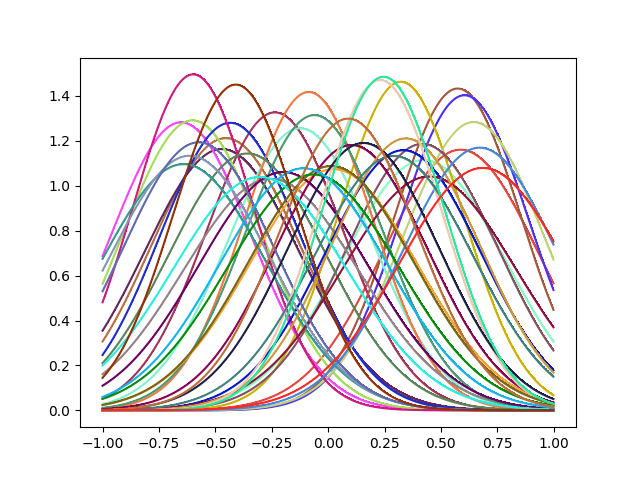}
	\caption{Initial opinion distribution}
	\label{fig: ic opinion distr}
	\medskip
	\small
	In this picture are represented the opinion distributions at initial time of the 40 agents considered for the simulation. Each distribution is described by a truncated Gaussian function, mean and variance of the Gaussian functions are independently uniform random distributed respectively in the intervals $[-0.7,0.7]$ and $[0.07, 0.15]$.
\end{figure}

\begin{figure}
	\centering
	\subfloat[$\rho_\textit{loc}=3.5$]{\includegraphics[width = 0.33\textwidth]{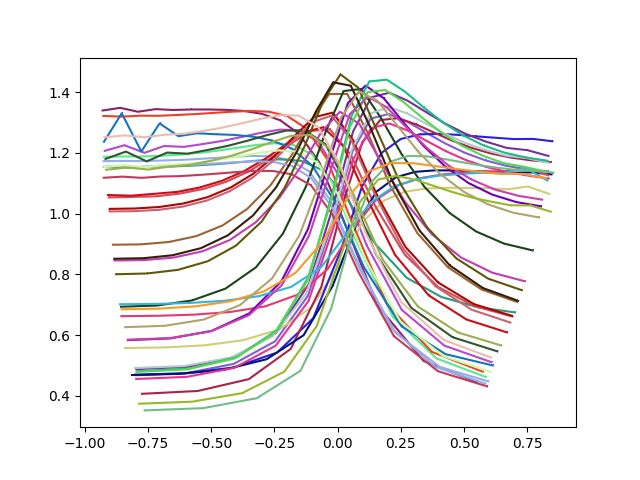}}
	\subfloat[$\rho_\textit{loc}=5$]{\includegraphics[width = 0.33\textwidth]{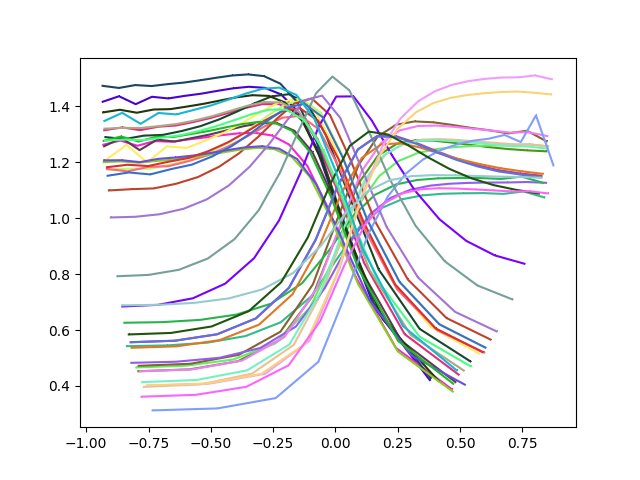}}
	\subfloat[$\rho_\textit{loc}=10$]{\includegraphics[width = 0.33\textwidth]{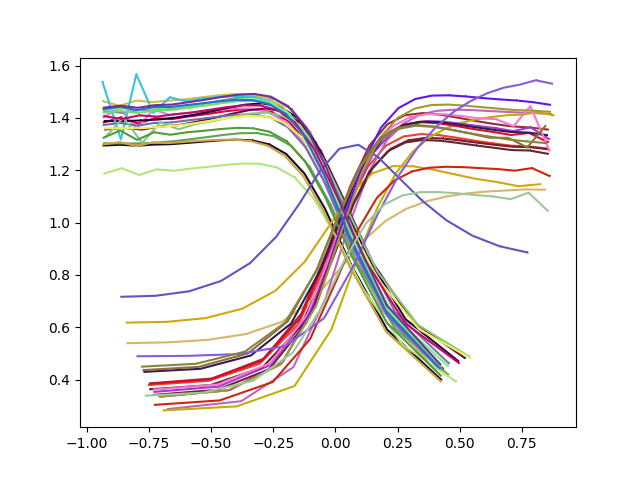}}
	\caption{Final opinion distribution}
	\label{fig: isto 1}
	\medskip
	\small
	We observe how the distributions are more and more concentrated either on the positive or negative side as the radius increases. Due to the diffusion the distributions tend to flatten once that they are concentrated on one of the two sides.\\
	The sharp oscillations close to the extreme values are due to the low resolution of the numerical partition of $\Omega$.
\end{figure}

Despite the common sense, we observe that in this case a wider range of interaction does not bring the system to a global consensus but to a more radicalized distribution of the opinions. Given the initial opinion distribution as in Figure \ref{fig: ic opinion distr},
we can observe that the distribution evolves towards a more and more radicalized society as the radius increases, see Figure \ref{fig: isto 1}.

\begin{figure}
	\centering
	\subfloat[$\rho_\textit{loc}=3.5$]{\includegraphics[width = 0.33\textwidth]{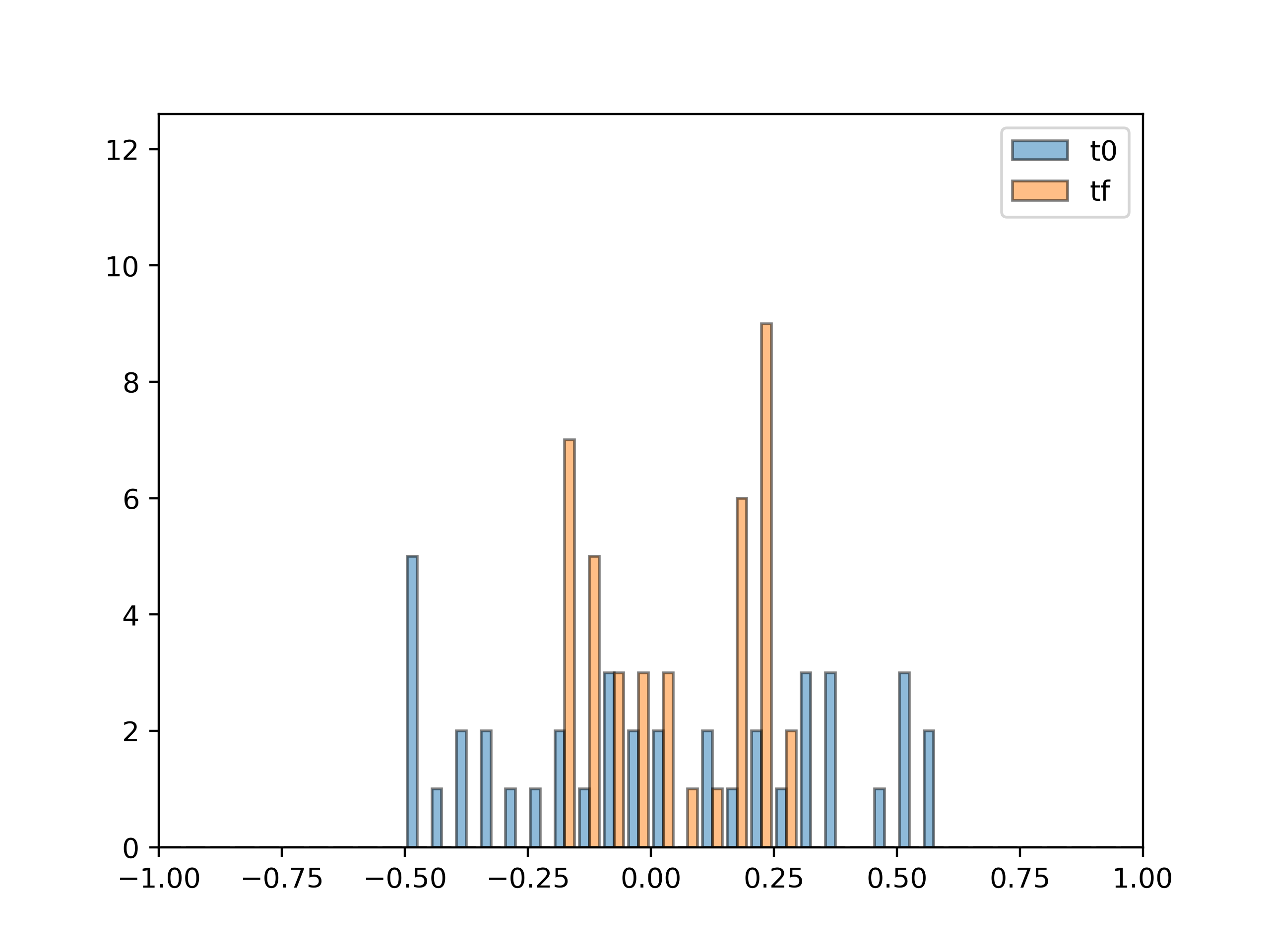}}
	\subfloat[$\rho_\textit{loc}=5$]{\includegraphics[width = 0.33\textwidth]{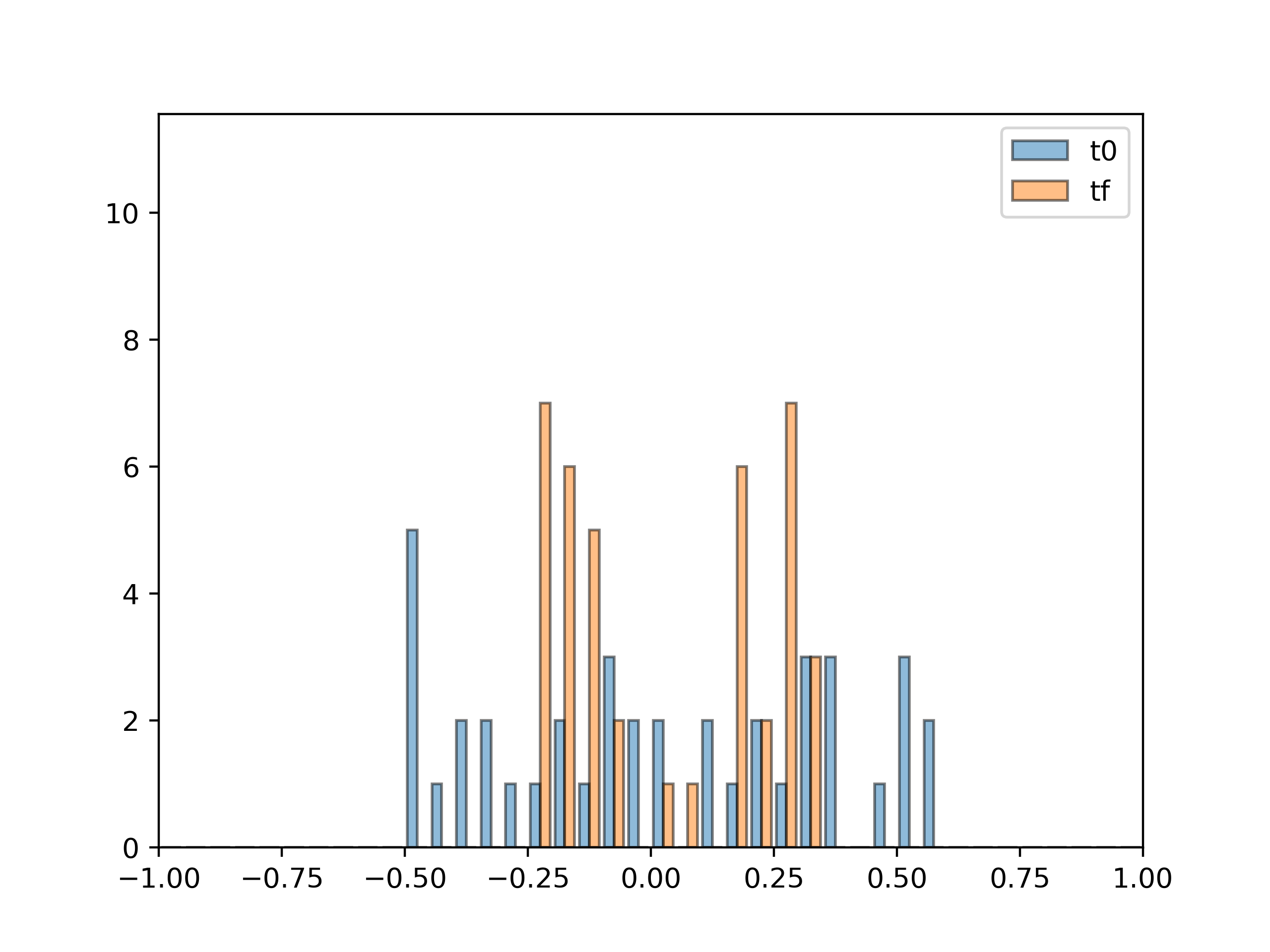}}
	\subfloat[$\rho_\textit{loc}=10$]{\includegraphics[width = 0.33\textwidth]{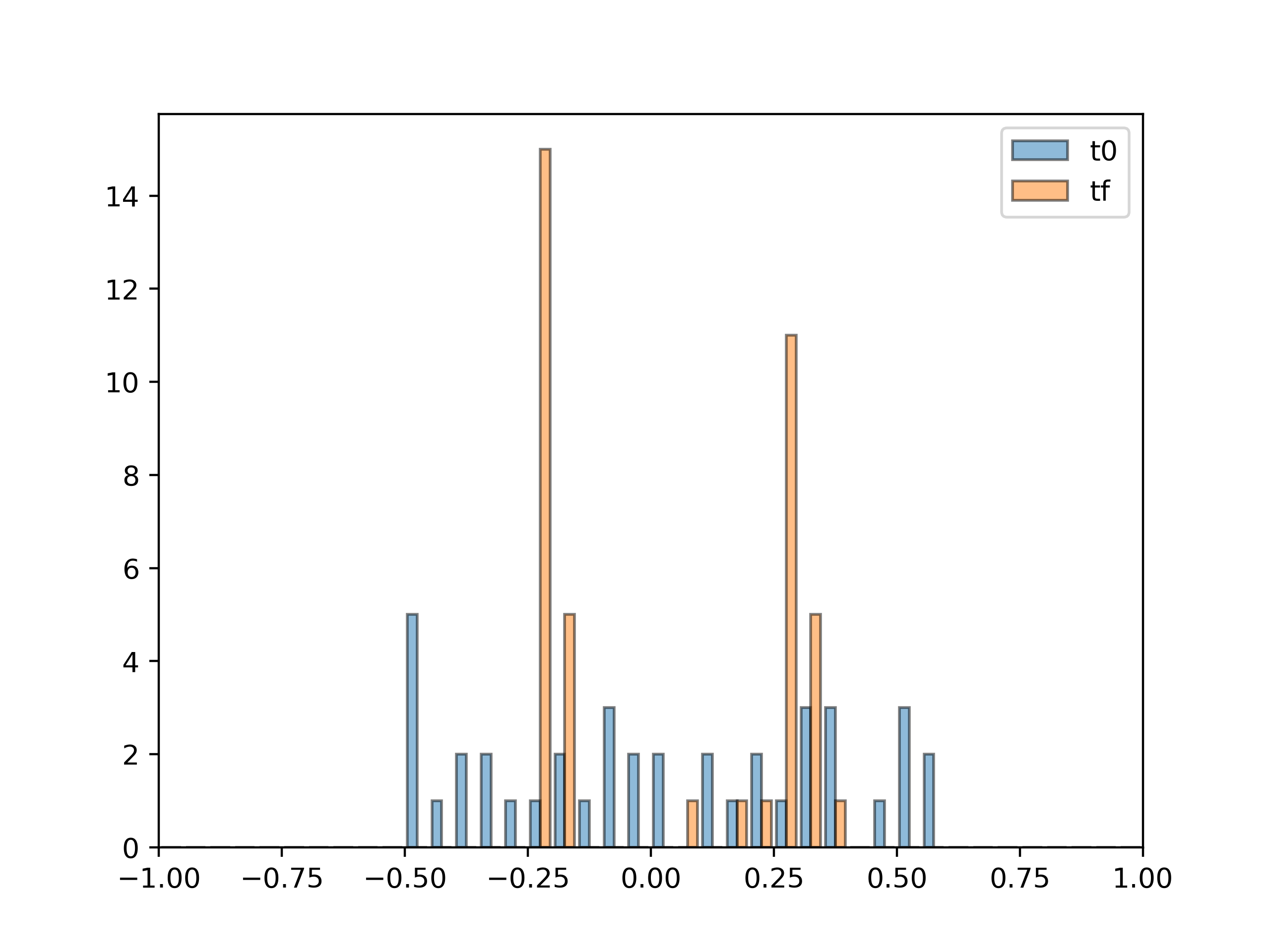}}
	\caption{Initial and final mean opinions distribution}
	\label{fig: isto 2}
	\medskip
	\small
	In this figure, in blue the distribution of the mean opinions at \color{blue}time $t=0$\color{black}, and in orange the mean opinions' distribution at \color{orange}final time \color{black} (which corresponds to the time showing a quasi stable status of the simulation result). We observe that the final distribution tends to have two peaks, which means that the opinions of the population are more and more split into two opinion's groups. However, they are also more close to the center. This means that we observe a sort of fragmentation and radicalization, but there is no polarization. 
\end{figure}

The same phenomenon is described by the histogram of the mean opinions distribution in Figure \ref{fig: isto 2}.

\begin{figure}
	\centering
	\subfloat[\color{blue}Blue]{\includegraphics[width = 0.33\textwidth]{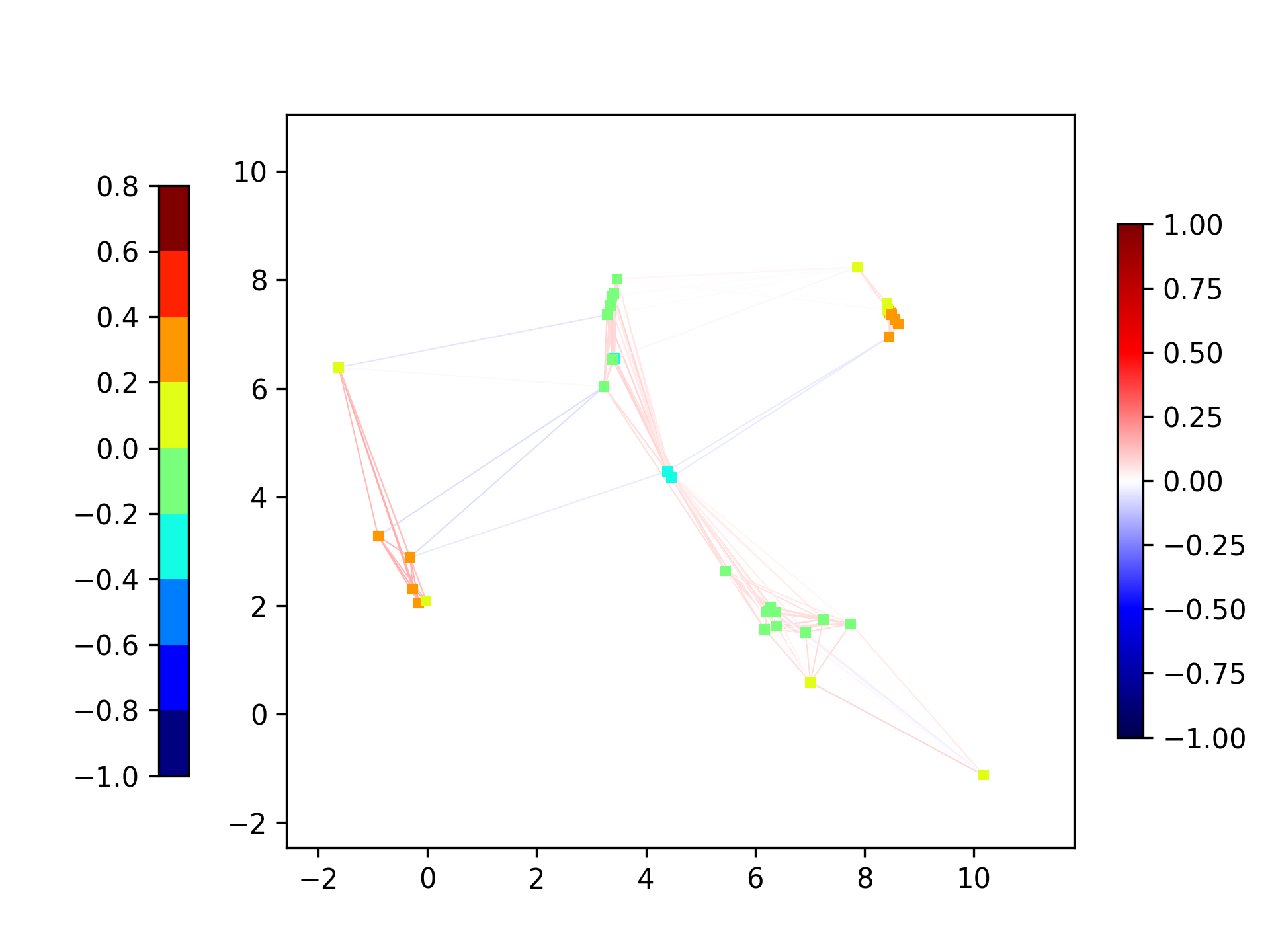}}
	\subfloat[Black]{\includegraphics[width = 0.33\textwidth]{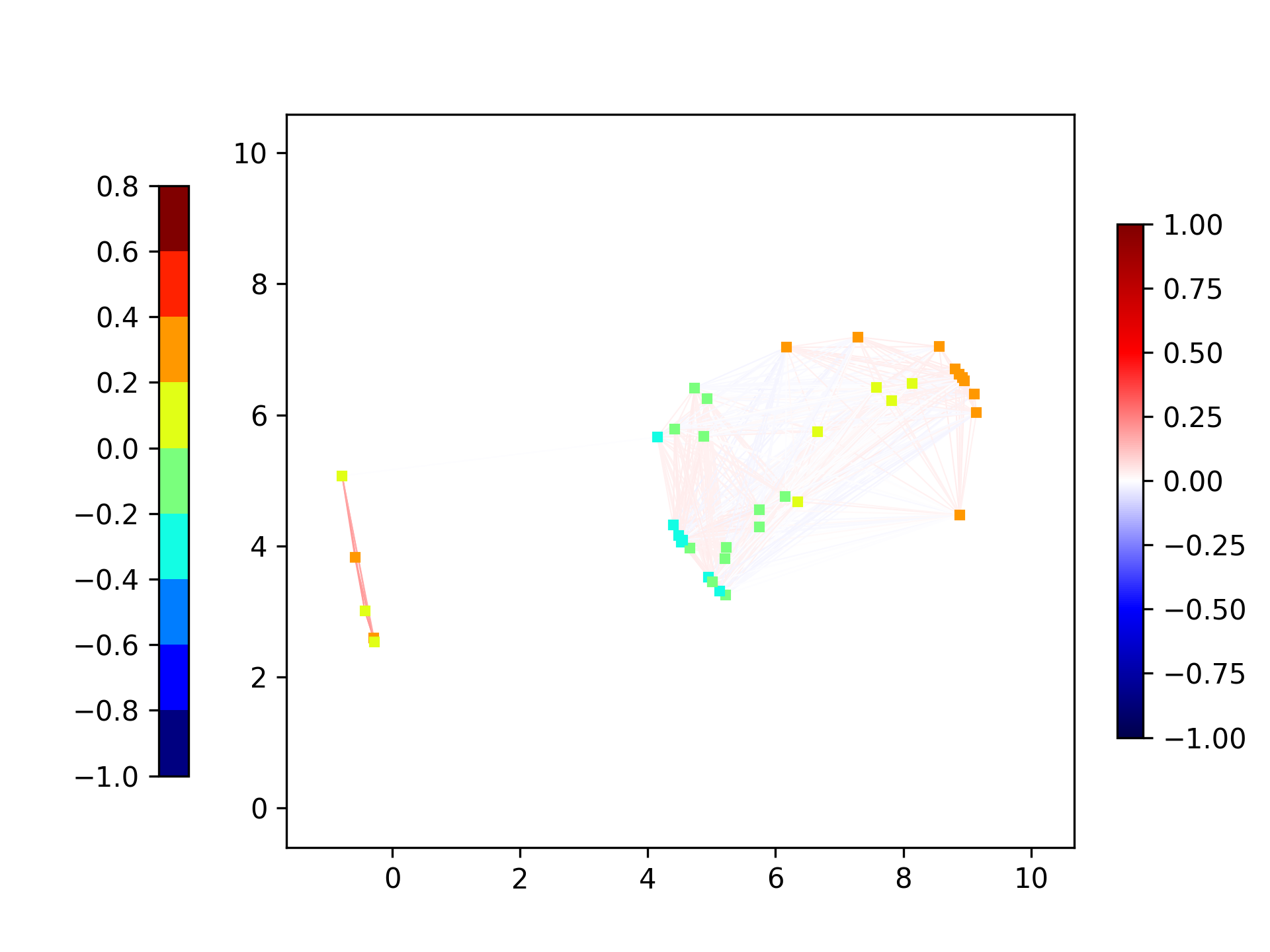}}
	\subfloat[\color{red}Red]{\includegraphics[width = 0.33\textwidth]{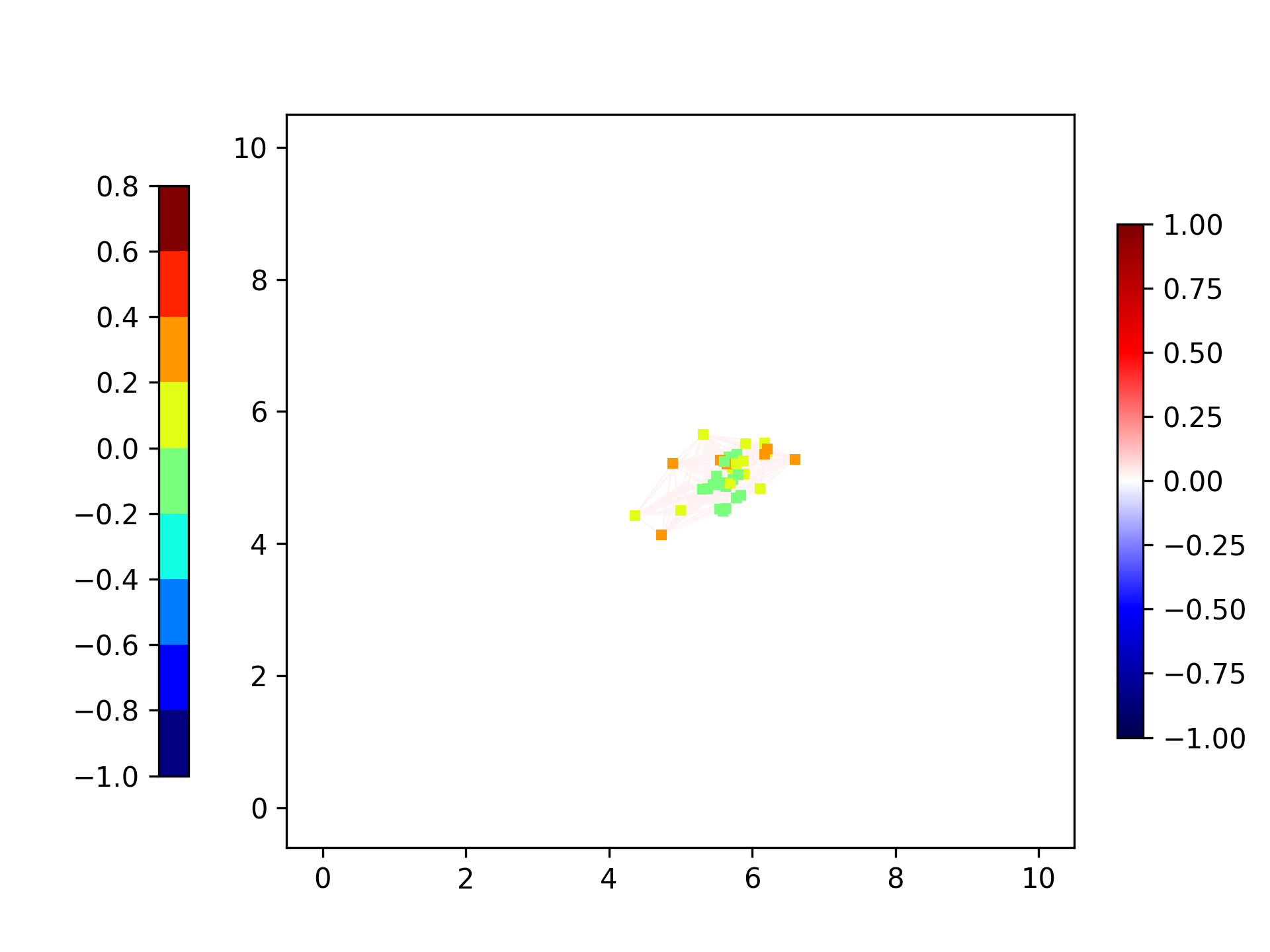}}
	\caption{Final network distribution}
	\label{fig: network}
	\medskip
	\small
	The olive function has not been plotted because it describes an extreme behaviour, all the agents collapse very fast into a unique point.\\ \color{blue}Blue\color{black}: $r_\textit{f}=0.15$, $r_\textit{a}=0.20$, $r_\textit{r}=0.30$, $r_\textit{l}=0.40$.\\ Black: $r_\textit{f}=0.25$, $r_\textit{a}=0.34$, $r_\textit{r}=0.36$, $r_\textit{l}=0.65$.\\ \color{red} Red\color{black}: $r_\textit{f}=0.30$, $r_\textit{a}=0.45$, $r_\textit{r}=0.55$, $r_\textit{l}=0.70$.
\end{figure}

It is interesting to compare the final network varying the parameters of the attitude areas. We fix the interaction radius $\rho_\textit{loc} = 5$ and we compare the first three functions plotted in Figure \ref{fig: scala}. In this case we observe that a more open minded society moves to a less fragmented network, see Figure \ref{fig: network}. 

\subsection{Polarization and fragmentation} 
While considering only the operator $\mathbf{K}^\mij$, i.e.~setting $\mathbf{A}^\mi=0$, we observe the phenomenon of the polarization. 

\begin{figure}
	\centering
	\subfloat[$\rho_\textit{loc}=3.5$]{\includegraphics[width = 0.33\textwidth]{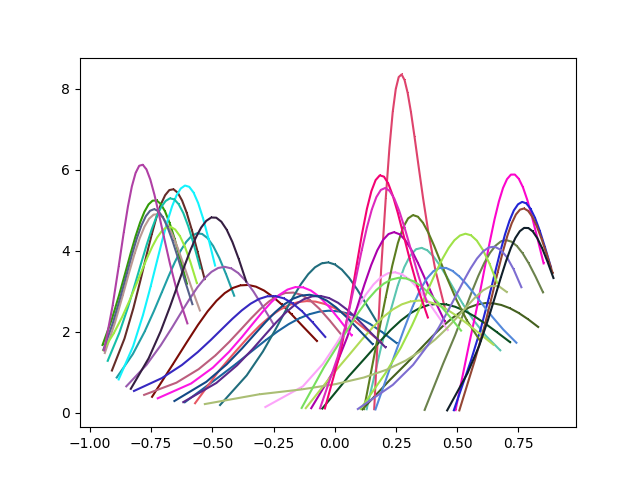}} 
	\subfloat[$\rho_\textit{loc}=5$]{\includegraphics[width = 0.33\textwidth]{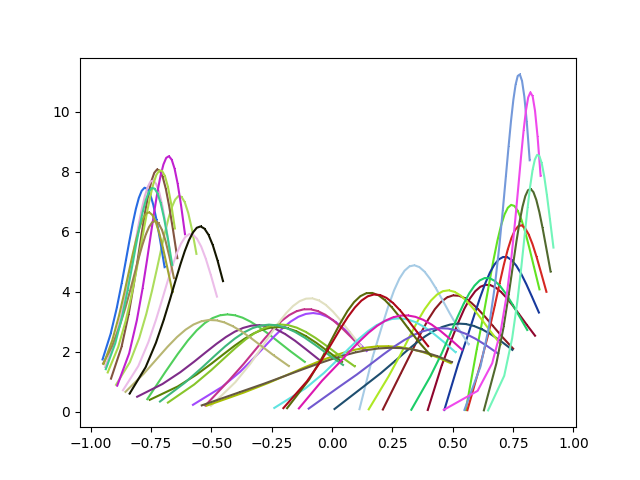}}
	\subfloat[$\rho_\textit{loc}=10$]{\includegraphics[width = 0.33\textwidth]{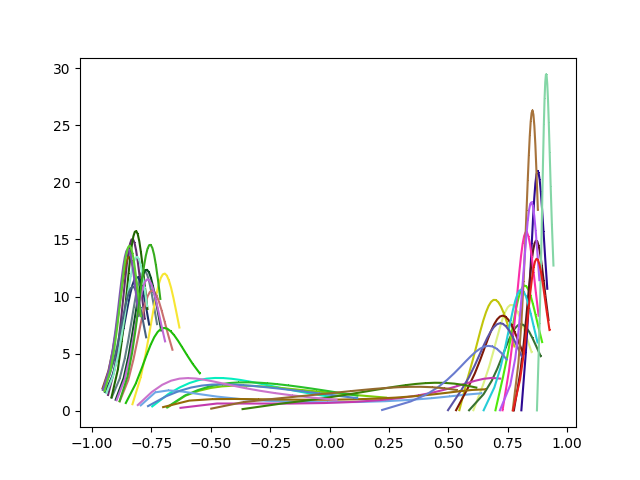}}
	\caption{Polarization while increasing the radius of interaction. Attitude function \color{blue}Blue\color{black}: $r_\textit{f}=0.15$, $r_\textit{a}=0.20$, $r_\textit{r}=0.30$, $r_\textit{l}=0.40$.}
	\label{fig: cresce raggio 1}
	\medskip
	\small
\end{figure}

It is interesting to observe how the connectivity of the population is the real driver of the polarization, instead of the \textit{open-mindedness}. We now compare the results of the simulation with fixed open-mindedness. We chose the most close minded population, i.e.~the one described by the blue function in Figure \ref{fig: scala}. Increasing the radius of the interaction, the opinion gets more and more polarized. In Figure \ref{fig: cresce raggio 1} we notice that a more connected society tends to cluster into extreme opinions. While a less connected society keeps a more sparse distribution of opinions. This is a not expected behaviour, usually a large or global interaction is related to a higher consensus. However, we introduced the operator $\mathbf{K}^\mij$ as the one describing the interaction on the social network, and so this results fits with the dynamics that we can observe nowadays considering the topics discussed mainly on the social platforms.

\begin{figure}
	\centering
	\subfloat[$\rho_\textit{loc}=3.5$]{\includegraphics[width = 0.33\textwidth]{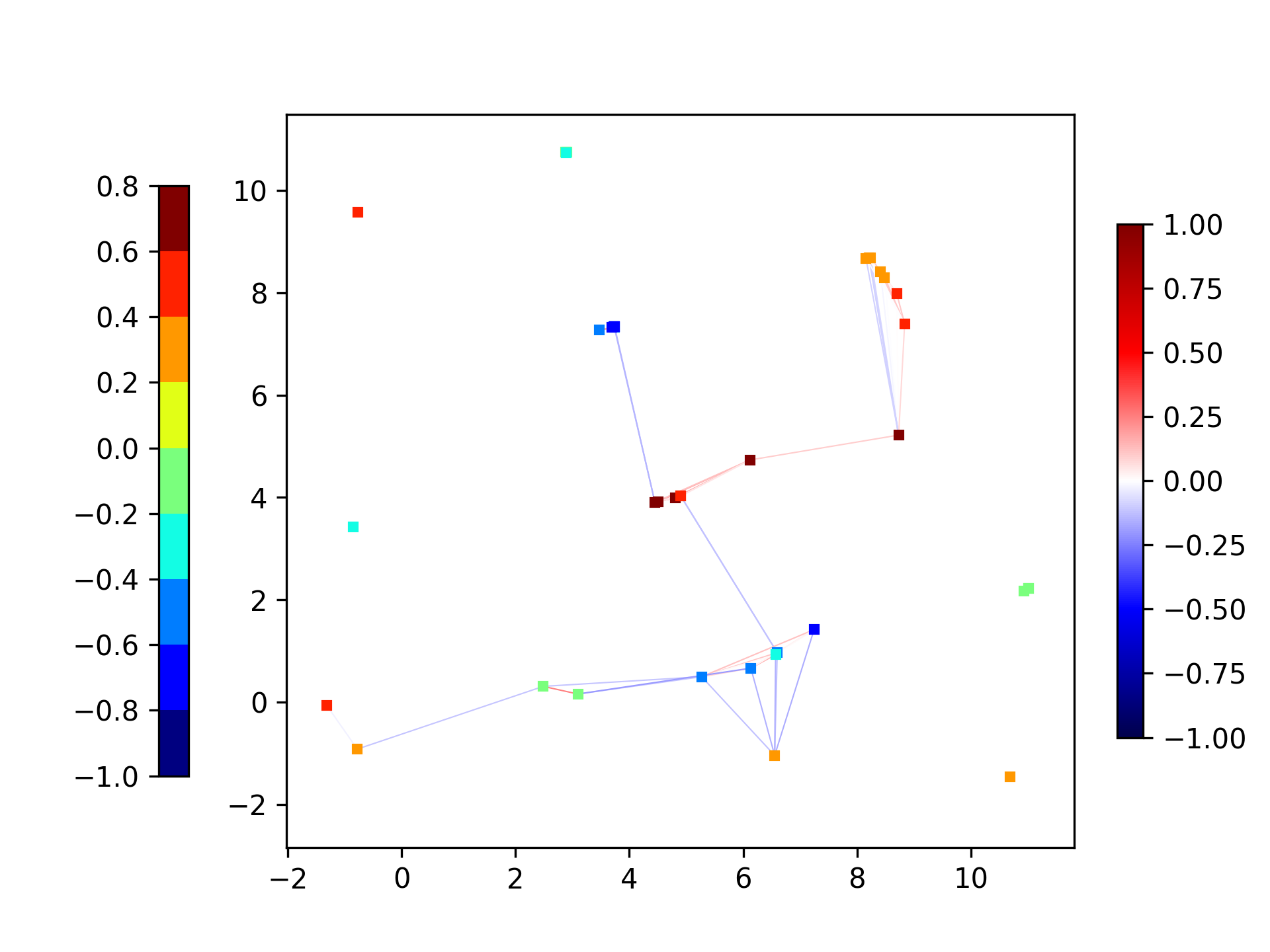}}
	\subfloat[$\rho_\textit{loc}=5$]{\includegraphics[width = 0.33\textwidth]{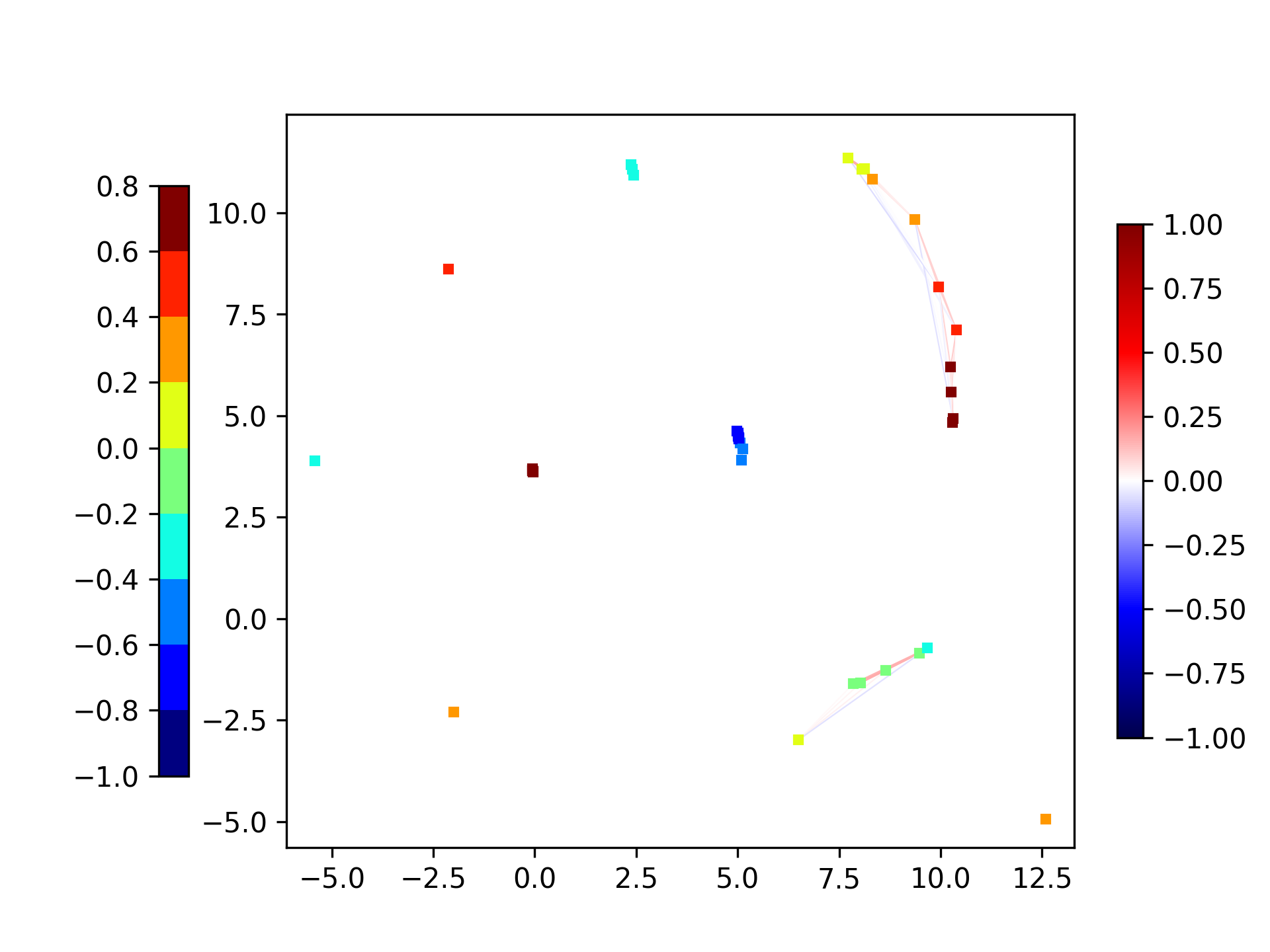}}
	\subfloat[$\rho_\textit{loc}=10$]{\includegraphics[width = 0.33\textwidth]{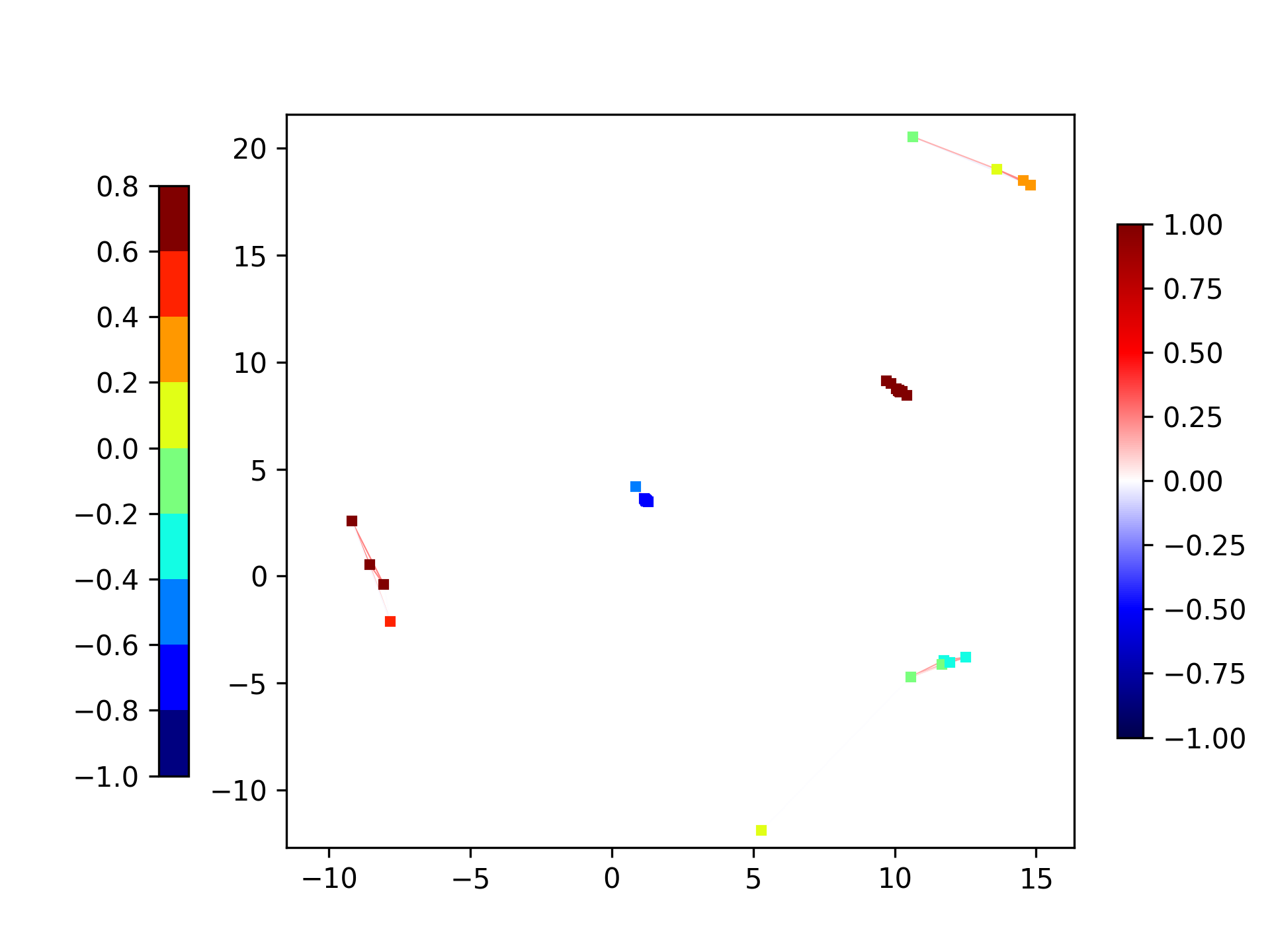}}
	\caption{Network fragmentation and opinion homogeneity}
	\label{fig: cresce raggio 2}
	\medskip
	\small
	Network while increasing the radius of interaction and keeping the same attitude function, i.e.~\color{blue}Blue\color{black}: $r_\textit{f}=0.15$, $r_\textit{a}=0.20$, $r_\textit{r}=0.30$, $r_\textit	{l}=0.40$.
\end{figure}

If we observe the evolution of the network, it seems that increasing the radius of interaction does not really affect the fragmentation of the population, but it plays a role in the opinion homogeneity of the network clusters. In Figure \ref{fig: cresce raggio 2}, at same fixed time, a larger radius of interaction brings to a wider network space, and the groups of connected agents show a stronger homogeneity of the opinion.

\subsection{Conclusions, interpretations, and possible follow up}
The goal of this model is to introduce the study of the processes describing the interaction on social networks and social media. We mainly focused on the role of the attitude areas and on that of the radius of interaction. The dynamics ruling the opinion formation of agents interacting on social platforms is different from the one described by models based on \textit{alignment}, \textit{averaged consensus}, or Cucker-Smale with positive communication rate.\\ Recently, sociologists and philosophers described the epistemic processes of the hyper connected society typical of the last two decades. The high amount of interactions and notions, together with their high frequency, modified the way how we create and reinforce our beliefs. Authors like Nguyen, see \cite{nguyen_2020}, explain how the network and the opinion distance are the discriminant for different epistemic processes. In our model we describe these two aspects through the definition of the attitude areas and through the dynamic of the network, which takes into account the distance on the network and the distance of the opinions. 

All together, the results obtained with the simulations show that the attitude area approach and the euclidean network structure describe a behaviour similar to the one observed on the social platforms. In particular, we observe the polarization typical of the social networks, and the fragmentation of the population in clusters with a strong opinion homogeneity.

Weakening the assumptions on the diffusion mobility $\mathbf{A}^\mi$ would make possible to model more precisely the role of social media. In particular, asking $\mathbf{A}^\mi$ to have Lipschitz first derivative is a very strong assumption that does not allow a realistic description of the interaction. Another step towards a more applied direction concerns the \textit{multi agent} description of the population. We already consider different \textit{masses} - i.e.~$\sigma^\mi$ - but it would be possible to define also a continuous family of attitude areas functions and its distribution among the population.

\section*{Acknowledgments}

The research of SF is  supported by the Ministry of University and Research (MIUR), Italy under the grant PRIN 2020- Project N. 20204NT8W4, Nonlinear Evolutions PDEs, fluid dynamics and transport equations: theoretical foundations and applications. The research of SF and GF is supported by the Italian INdAM project N. E55F22000270001 ``Fenomeni di trasporto in leggi di conservazione e loro applicazioni''. SF is also supported by University of L'Aquila 2021 project 04ATE2021 - ``Mathematical Models For Social Innovations: Vehicular And Pedestrian Traffic, Opinion Formation And Seismology.''

\printbibliography

@article{T1,
	title={Kinetic models of opinion formation},
	author={Toscani, G.},
	journal={Comm. Math. Sci.},
	volume={4},
	number={3},
	pages={481--496},
	year={2006}
}

@article {russo1,
	author = {Russo, G.},
	title = {Deterministic diffusion of particles},
	journal = {Comm. on Pure and Applied Mathematics},
	volume = {43},
	pages = {697-733},
	year = {1990}}

@article{FaRa,
	author = {Fagioli, S. and Radici, E.},
	title = {Solutions to aggregation-diffusion equations with nonlinear mobility constructed via a deterministic particle approximation},
	journal = {Math. Mod. and Meth. in App. Sci.},
	volume = {28},
	number = {09},
	pages = {1801-1829},
	year = {2018}
}

@article{RoSa,
	AUTHOR = {Rossi, R. and Savar{\'e}, G.},
	TITLE = {Tightness, integral equicontinuity and compactness for
		evolution problems in {B}anach spaces},
	JOURNAL = {Ann. Sc. Norm. Super. Pisa Cl. Sci. (5)},
	FJOURNAL = {Annali della Scuola Normale Superiore di Pisa. Classe di
		Scienze. Serie V},
	VOLUME = {2},
	YEAR = {2003},
	NUMBER = {2},
	PAGES = {395--431},
	ISSN = {0391-173X},
	MRCLASS = {46E30 (28A20 35K55 46N20)},
	MRNUMBER = {2005609},
	MRREVIEWER = {Agnieszka Ka{\l}amajska},
}

@book {AGS,
	AUTHOR = {Ambrosio, L. and Gigli, N. and Savar{\'e}, G.},
	TITLE = {Gradient flows in metric spaces and in the space of
		probability measures},
	SERIES = {Lectures in Mathematics ETH Z\"urich},
	EDITION = {Second},
	PUBLISHER = {Birkh\"auser Verlag},
	ADDRESS = {Basel},
	YEAR = {2008}
}

@book {V1,
	AUTHOR = {Villani, C.},
	TITLE = {Topics in optimal transportation},
	SERIES = {Graduate Studies in Mathematics},
	VOLUME = {58},
	PUBLISHER = {American Mathematical Society},
	ADDRESS = {Providence, RI},
	YEAR = {2003}}

@book{S,
	author = {Santambrogio, F. },
	title = { Optimal Transport for Applied Mathematicians },
	series = { Progress in Nonlinear Differential Equations and Their Applications },
	volume = {86},
	PUBLISHER = {Birkh\"auser Verlag},
	ADDRESS = {Basel},
	YEAR = {2015}}

@article{CT,
title={Wasserstein metric and large--time asymptotics of nonlinear diffusion equations},
author={Carrillo, J. A. and Toscani, G.},
journal={New Trends in Mathematical Physics, (In Honour of the Salvatore Rionero 70th Birthday)},
pages={234--244},
year={2005}
}

@article{MoTa,
author = {Motsch, S. and Tadmor, E.},
title = {Heterophilious Dynamics Enhances Consensus},
journal = {SIAM Review},
volume = {56},
number = {4},
pages = {577-621},
year = {2014}}

@article{BoLo,
Author = {Borra, D. and Lorenzi, T.},
Journal = {Zeitschrift f{\"u}r angewandte Mathematik und Physik},
Number = {3},
Pages = {419--437},
Title = {A hybrid model for opinion formation},
Volume = {64},
Year = {2013}}

@book{BeMaTo,
author = {Bellomo, N. and Ajmone Marsan, G. and Tosin, A.}, 
title={Complex Systems and Society. Modeling and Simulation. SpringerBriefs in Mathematics}, 
publisher = {Springer},
year = {2013} }

@article{CaFoLo, 
author= {Castellano, C. and Fortunato, S. and Loreto, V.}, 
title = {Statistical physics of social dynamics.},
journal= {Review of Modern Physics}, 
volume = {81},
number={2},
pages={591-646}, 
year = {2009}}

@book{NaPaTo,
author = {Naldi, G. and  Pareschi, L. and Toscani, G.},
title = { Mathematical Modeling of Collective Behavior in Socio-Economic and Life Sciences},
PUBLISHER = {Birkh\"auser},
ADDRESS = {Boston},
YEAR = {2010}}

@book{PaTobook,
author = {Pareschi, L. and Toscani, G.},
title = { Interacting Multiagent Systems. Kinetic Equations and Monte Carlo Methods. },
PUBLISHER = {Oxford University Press},
YEAR = {2013}}

@book{Ga,
author = {Galam, S.},
title = { Sociophysics: a physicists modeling of psycho-political phenomena (understanding complex systems)},
PUBLISHER = {Springer},
YEAR = {2012}}

@article{St, 
author= {Strogatz, S. H}, 
title = {Exploring complex networks},
journal= {Nature}, 
volume = {410},
number={6825},
pages={268-276}, 
year = {2001}}

@article{YaRoSc, 
author= {Yardi, S. and Romero, D. and Schoenebeck, G.}, 
title = {Detecting spam in a Twitter
	network},
journal= {First Monday}, 
volume = {15},
number={1},
year = {2009}}

@article{KlShSh, 
author= {Klein, A. and Ahlf, H. and Sharma, V. }, 
title = { Social activity and structural centrality
	in online social networks},
journal= {Telematics and Informatics}, 
volume = {32},
number={2},
pages={321-332}, 
year = {2015}}

@article{Sz,
author= {Sznajd-Weron, K. and Sznajd, J.},
year =  {2000},
title =  {Opinion evolution in closed community.},
journal={ Int. J. Mod. Phys. C},
volume ={11},
pages={1157?1165}}

@article{SlLa,
author = {Slanina, F. and Lavi\v{c}ka, H.},
year= {2003},
title = {Analytical results for the Sznajd model of opinion formation.},
journal={ Eur.Phys. J. B},
volume= {35}, 
pages={279?288}}

@inbook{AlPaToZa,
title = {Recent advances in opinion modeling: control and social influence},
author = {Albi, G. and Pareschi, P. and Toscani, G. and Zanella, M.},                                                   
bookTitle = {Active Particles, Volulme 1. Advances in Theory, Models, and Applications}, 
editor = {Bellomo, N. and Degond, P. and Tadmor, E.},
publisher={Birkh\"auser - Springer},
pages={49-98},
year ={ 2017}}

@article{Be,
author = {Ben-Naim, E.},
title={Opinion dynamics: rise and fall of political parties}, journal={Europhysics Letters}, 
volume={69},
number={5},
pages={671},
year= {2005}}

@incollection{AlPaZa3,
author = {Albi, G. and Pareschi, L. and Zanella, M.},
title = {On the Optimal Control of Opinion Dynamics on Evolving Networks},
booktitle = {System Modeling and Optimization. CSMO 2015. IFIP Advances in Information and Communication Technology}, 
editors= {Bociu, L. and Désidéri, J.A. and Habbal, A.}, 
volume = {494},
pages = {58-67},
publisher = {Springer},
address= {Cham},
year = {2016}}

@article{LaMa,
title={Opinion propagation on social networks: a Mathematical Standpoint},
author={Lavenant, H. and Maury, B.},
journal = {preprint},
pages = {53},
year={2019}}

@article{DiFSt,
title = "Convergence of the follow-the-leader scheme for scalar conservation laws with space dependent flux",
journal = "Discrete \& Continuous Dynamical Systems - A",
volume = "40",
number = "1",
pages = "233-266",
year = "2020",
author = "Di Francesco, M. and Stivaletta, G."}

@article{DiFFaRa,
journal = {Journal of Differential Equations},
volume =  {266}, 
number = {5}, 
year = {2019}, 
pages =  {2830-2868},
title = {Deterministic particle approximation for nonlocal transport equations with nonlinear mobility},
author = {Di Francesco, M. and Fagioli, S. and Radici, E.}}

@article{DiFRo,
author = {Di Francesco, M. and Rosini, M.D.},
title = {Rigorous derivation of nonlinear scalar conservation laws from follow-the-leader type models via many particle limit}, journal = {Archive for rational mechanics and analysis}, 
volue = {217}, 
issue = {3}, 
pages = {831-871},
year = {2015}}

@article{Burger_Net,
author = {Burger, M.},
title = {Kinetic equations for processes on co-evolving networks},
journal ={Kinetic and Related Models}, 
volume = {15},
number = {2}, 
pages = {187-212}, 
year = {2022}}

@article{Burger_OF,
author = {Burger, M.},
title = {Network structured kinetic models of social interactions},
journal ={Vietnam J Math.}, 
volume = {49},
number = {3}, 
pages = {937-956}, 
year = {2021}}

@article{FagRad,
author = {Fagioli, S. and Radici, E.},
title = {Opinion formation systems via deterministic particles approximation},
journal ={Kinetic and Related Models}, 
volume = {14},
number = {1}, 
pages = {45-76}, 
year = {2020}}

@article{FagTse,
author = {Fagioli, S. and Tse, O.},
title = {On gradient flow and entropy solutions for nonlocal transport equations with nonlinear mobility},
journal ={Nonlinear Analysis}, 
volume = {221},
pages = {112904}, 
year = {2022}}

@article {diFFRo2017,
AUTHOR = {Di Francesco, M. and Fagioli, S. and Rosini, M. D.},
TITLE = {Deterministic particle approximation of scalar conservation
	laws},
JOURNAL = {Boll. Unione Mat. Ital.},
FJOURNAL = {Bollettino dell'Unione Matematica Italiana},
VOLUME = {10},
YEAR = {2017},
NUMBER = {3},
PAGES = {487--501},
ISSN = {1972-6724},
MRCLASS = {35L65 (90B20)},
MRNUMBER = {3691810}
}

@article{GT1,
Author = {L. Gosse and G. Toscani},
Journal = {SIAM J. Numer. Anal.},
Title = {Identification of asymptotic decay to self-similarity for one- dimensional filtration equations},
Volume = {43},
Year = {2006},
Pages = {2590–2606}}

@article{BHJ11,
title = {A moving mesh finite element algorithm for the adaptive solution of time-dependent partial differential equations with moving boundaries},
journal = {Journal of Computational and Applied Mathematics},
volume = {54},
pages = {450-469},
year = {2004},
author = {Baines, M.J. and Hubbard, M.E. and Jimack, P.K.}}

@article{BHR96,
title = {Moving mesh methods for problems with blow-up},
journal = {SIAM Journal on Scientific Computing},
volume = {17},
pages = {305-327},
year = {1996},
author = {Budd, C. and Huang, W. and Russell, R.}}

@article{BHR09,
title = {Adaptivity with moving grids},
journal = {Acta Numerica},
volume = {18},
pages = {111-241},
year = {2009},
author = {Budd, C. and Huang, W. and Russell, R.}}

@article{SMR01,
title = {A moving mesh method for one-dimensional hyperbolic conservation laws},
journal = {SIAM Journal on Scientific Computing},
volume = {22},
pages = {1791-1813},
year = {2001},
author = {Stockie, J. and Mackenzie, J. and Russell, R.}}

@article{CHR02,
title = {A moving mesh method based on the geometric conservation law},
journal = {SIAM Journal on Scientific Computing},
volume = {24},
pages = {118-142},
year = {2002},
author = {Cao, W. and Huang, W. and Russell, R.}}

@article{CHR03,
title = {Approaches for generating moving adaptive meshes: location versus velocity},
journal = {Applied Numerical Mathematics},
volume = {47},
number = {2},
pages = {121-138},
year = {2003},
author = {Cao, W. and Huang, W. and Russell, R.}}

@article{begby_2022, title={From Belief Polarization to Echo Chambers: A Rationalizing Account}, journal={Episteme}, publisher={Cambridge University Press}, author={E. Begby}, year={2022}, pages={1–21}}

@article{nguyen_2020, 
title={Echo chambers and epistemic bubbles}, 
volume={17}, 
number={2}, 
journal={Episteme}, 
publisher={Cambridge University Press},
author={Nguyen, C. Thi}, 
year={2020}, 
pages={141–161}}

@book{ThHaKl, 
title={Introduction to the Theory of Complex Systems},  
publisher={Oxford University Press},
author={S.Thurner, S. and Hanel, R. and Klimek, P.}, 
year={2018}}

@article{benatti2020da,
author = "Benatti, A. and de Arruda, H. F. and Silva, F. N. and Comin, C. H.",
title = "da Fontoura Costa",
journal = "L. (2020). {Opinion} diversity and social bubbles in adaptive {Sznajd} networks. {Journal} of {Statistical} {Mechanics}: {Theory} and {Experiment}, 2020(2)",
address = "Opinion diversity and social bubbles in adaptive Sznajd networks. Journal of Statistical Mechanics",
publisher = "Theory and Experiment",
pages = "023407",
year = "2020"
}

@incollection{nigam2018onem,
author = "A. Nigam and K. Shin and A. Bahulkar and B. Hooi and D. Hachen and B. K. Szymanski and C. Faloutsolos and N. V. Chawla",
title = "ONE-{M}: modeling the co-evolution of opinions and network connections",
booktitle = "Joint European Conference on Machine Learning and Knowledge Discovery in Databases (2018)",
pages = "122--140",
year = "2018"
}

@article{tur2018coevolution,
author = "E. M. Tur and J. M. Azagra-Caro",
title = "The coevolution of endogenous knowledge networks and knowledge creation",
journal = "Journal of Economic Behavior and Organization",
volume = "145.",
pages = "424--434",
year = "2018"
}

@incollection{kohne2020role,
author = "J. Kohne and N. Gallagher and Z. M. Kirgil and R. Paolillo and L. Padmos and F. Karimi",
title = "The role of net- work structure and initial group norm distributions in norm conflict",
booktitle = "Computational Conflict Research, Springer, Cham (2020)",
pages = "113--140",
year = "2020"
}

@incollection{albi2017continuum,
author = {Albi, G. and Burger, M. and Haskovec, J. and Markowich, P. and Schlottbom, M.},
title = {Continuum modeling of biological network formation},
booktitle = {System Modeling and Optimization. CSMO 2015. IFIP Advances in Information and Communication Technology}, 
volume = {494},
pages = {58-67},
publisher = {Springer},
address= {Cham},
year = {2016}}

@article{coppini2020law,
author = "F. Coppini and H. Dietert and G. Giacomin",
title = "A law of large numbers and large deviations for interacting diffusions on ErdoesRenyi graphs",
journal = "Stochastics and Dynamics 20 (2020)",
pages = "2050010.",
year = 2020
}

@article{delattre2016lucon,
author = "S. Delattre and G. Giacomin and E. Lucon",
title = "Lucon {A} note on dynamical models on random graphs and {FokkerPlanck} equations",
journal = "Journal of Statistical Physics",
volume = "165",
pages = "785--798",
year = 2016
}

@article{Lackey_of,
author = {J. Lackey},
isbn = {9780198863977},
title = "{Echo Chambers, Fake News, and Social Epistemology}",
booktitle = "{The Epistemology of Fake News}",
publisher = {Oxford University Press},
year = {2021},
month = {06}}

@incollection{Bernecker_of,
author = {S. Bernecker},
isbn = {9780198863977},
title = "{An Epistemic Defense of News Abstinence}",
booktitle = "{The Epistemology of Fake News}",
publisher = {Oxford University Press},
year = {2021},
month = {06}}

@article{O_Connor_of,
title={Endogenous epistemic factionalization},
author={J.O. Weatherall and C. O’Connor},
journal={Synthese},
volume={198},
number={Suppl 25},
pages={6179--6200},
year={2021},
publisher={Springer}}

@article{Rosenstock_of, 
title={In Epistemic Networks, Is Less Really More?}, 
volume={84}, 
number={2}, 
journal={Philosophy of Science}, 
publisher={Cambridge University Press}, 
author={S. Rosenstock and J. Bruner and C. O’Connor}, 
year={2017}, 
pages={234–252}}

@article{wolfram_of,
title={On evolving network models and their influence on opinion formation},
author={Nugent, A. J. and Gomes, S. N. and Wolfram, M.-T.},
journal={arXiv preprint arXiv:2305.09483},
year={2023}
}

@book{bouchut2000series,
author = "F. Bouchut and F. Golse and M. Pulvirenti: Kinetic Equations and Asymptotic Theory",
title = "Series in Applied Mathematics, 4, Gauthier-Villars, Paris",
year = 2000
}

@book{cercignani1969mathematical,
author = "C. Cercignani",
title = "Mathematical Methods in Kinetic Theory",
publisher = "Plenum Press, New York",
year = 1969
}

@article{burger2014partial,
author = "M. Burger and L. Caffarelli and P. A. Markowich",
title = "Partial differential equation models in the socio- economic sciences",
journal = "Phil. Trans. Royal Society",
volume = "372",
pages = "20130406.",
year = 2014
}

@book{naldi2010eds,
author = "G. Naldi and L. Pareschi and G. Toscani",
title = "eds., Mathematical Modeling of Collective Behavior in Socio-Economic and Life Sciences",
publisher = "Springer, New York",
year = 2010
}

@article{pareschi2014wealth,
author = "L. Pareschi and G. Toscani",
title = "Wealth distribution and collective knowledge: a {Boltzmann} approach",
journal = "Philosophical Transactions of the Royal Society A: Mathematical, Physical and Engineering Sciences",
number = "372",
pages = "20130396",
year = 2014
}

@article{boud2009,
author = "Boudin, L. and Salvarani, F.",
title = "A kinetic approach to the study of opinion formation.",
journal = "ESAIM: Mathematical Modelling and Numerical Analysis",
number = "43",
pages = "507-522",
year = 2009
}

\end{document}